\documentclass{article}

\usepackage{arxiv}

\usepackage[utf8]{inputenc} % allow utf-8 input
\usepackage[T1]{fontenc}    % use 8-bit T1 fonts
\usepackage{hyperref}       % hyperlinks
\usepackage{url}            % simple URL typesetting
\usepackage{booktabs}       % professional-quality tables
\usepackage{amsfonts}       % blackboard math symbols
\usepackage{nicefrac}       % compact symbols for 1/2, etc.
\usepackage{microtype}      % microtypography
\usepackage{lipsum}
\usepackage{fancyhdr}       % header
\usepackage{graphicx}       % graphics
\usepackage{lmodern}
\usepackage[style=numeric]{biblatex}
\usepackage{doi}
\usepackage{geometry}

\usepackage{subcaption}
\usepackage{enumerate}
\usepackage[shortlabels]{enumitem}
\usepackage{amssymb}
\usepackage[normalem]{ulem}
\usepackage{titlesec,bbm}
\usepackage{spverbatim} 
\usepackage{amsmath}
\usepackage{amsthm}
\usepackage{float}
\usepackage{multirow}
\usepackage{authblk}
\usepackage{algorithmic}
\usepackage{listings} 
\usepackage{tikz}
\usepackage{relsize}
\usepackage[title]{appendix}
\usepackage{bm}
\newcommand{\R}{\mathbb R}

\newcommand{\btheta}{\bm{\theta}}

\newcommand{\bfK}{\bm{K}}
\newcommand{\bfQ}{\bm{Q}}
\DeclareMathOperator{\sech}{sech}

\title{Transformed Physics-Informed Neural Networks for The Convection-Diffusion Equation} %%%%%%%%%%%%
\author{Jiajing Guan}
\author{Howard Elman}
\affil{University of Maryland}
\date{\today}
\begin{document}
+\maketitle

\begin{abstract}
Singularly perturbed problems are known to have solutions with steep boundary layers that are hard to resolve numerically. Traditional numerical methods, such as Finite Difference Methods (FDMs), require a refined mesh to obtain stable and accurate solutions. As Physics-Informed Neural Networks (PINNs) have been shown to successfully approximate solutions to differential equations from various fields, it is natural to examine their performance on singularly perturbed problems. The convection-diffusion equation is a representative example of such a class of problems, and we consider the use of PINNs to produce numerical solutions of this equation. We study two ways to use PINNS: as a method for correcting oscillatory discrete solutions obtained using FDMs, and as a method for modifying reduced solutions of unperturbed problems. For both methods, we also examine the use of input transformation to enhance accuracy, and we explain the behavior of input transformations analytically, with the help of neural tangent kernels.

% After identifying the inherent issues with PINNs, we explain their behaviors through the lens of neural tangent kernel. Finally, we introduce a simple but effective input transformation scheme, targeting the inability of PINNs in accurate approximating solutions to convection-diffusion equations. 

\end{abstract} %%%%%%%%%

% \bigskip

% \noindent \lipsum[1] \cite{1}

% $\,$

% $\,$

\section{Introduction}
% introduction on convection diffusion equation (singularly perturbed problem, boundary layer) and numerical methods

% On the other hand, these convection-diffusion equations are part of the singularly perturbed problem family. 

Singularly perturbed problems are differential equations that depend on a small positive parameter $\epsilon$ and whose solutions (or their derivatives) at the boundary contain a steep layer as $\epsilon$ approaches zero. The parameter $\epsilon$ is called the perturbation parameter \cite{roos2008robust}. Due to the steep boundary layer, stable and accurate numerical solutions of such equations are hard to obtain. Special numerical techniques have been developed to resolve the layers \cite{bakhvalov1969optimization,gartland1988graded,miller1996fitted,linss2003layer}. Adaptive meshes are commonly used to generate meshes that have high resolution near the boundary layer \cite{linss2003layer}. However, due to the fine resolution, computational costs may be high and there is an interest in find numerical techniques that can solve the equations stably and accurately without large computational costs. 

The convection-diffusion equation is a representative example in the class of singularly perturbed differential equations. Convection-diffusion equations are commonly used to describe the combined effects of two different physical processes: the diffusion of a content within a fluid, and the swift movement of the fluid that convects the content downstream. Such equations are used to model fluid flows in a wide range of disciplines including climate studies, biological systems, chemical processes, energy and astrophysics. Common numerical schemes for solving the convection-diffusion equation are finite difference methods (FDM). finite element methods (FEM) and finite volume methods (FVM). Nonetheless, without delicate numerical techniques or adaptive meshes, nonphysical oscillations will occur near the boundary layer. 

% introduction on neural network and PINN
In recent decades, due to developments in computing hardware, deep neural networks have gained attention as large models become feasible, and they have been advocated as suitable surrogate functions for PDEs. Physics-informed neural networks (PINNs) were introduced in \cite{RAISSI2019686}, approximating the solutions to a range of problems in computational science and engineering. They have been shown to be effective for models of fluid mechanics \cite{raissi2020hidden,sun2020surrogate}, bio-engineering \cite{sahli2020physics,kissas2020machine}, meta-material design \cite{fang2019deep,liu2019multi}, free boundary problems \cite{wang2021deep}, and many more fields. The benefit of using PINNs as surrogate functions is that, once trained, approximate solutions can be generated quickly with low cost. However, PINNs
using fully connected architectures may fail to be trained effectively, which can result in inaccurate surrogate approximations. We would like to explore if PINNs can avoid the difficulties associated with steep boundary layers and avoid nonphysical oscillations near boundary layers. We will also develop ways to make the training process robust. We are interested in exploring using PINNs as correctors for solutions of convection-diffusion equation obtained various ways, primarily focusing on one-dimensional convection-diffusion equations as test cases. 

% Thus, it is desirable to apply PINNs on convection-diffusion equations to investigate if PINN can avoid the issue with stiff boundary layer.\\

% introduction on structure UPDATE STRUCTURE
The outline of the paper is as follows. We lay the foundations of neural networks and PINNs in Section \ref{sec:nn}, where we also introduce the input transformation technique. We then review convection-diffusion equations in Section \ref{sec:cd}. We explore the use of PINNs to correct unphysical oscillations in discrete solutions obtained from FDM in Section \ref{sec:correcting_solution}, with and without input transformation. We investigate how PINNs perform in correcting reduced solutions and demonstrate the power of input transformations in Section \ref{sec:correcting_reduced}. We examine the effect of input transformations analytically in Section \ref{sec:input_transformation}. Finally, we employ similar techniques in correcting reduced solutions on a two-dimensional convection-diffusion test problem in Section \ref{sec:2dcd}.

\section{Neural Networks and Physics-Informed Neural Networks}
\label{sec:nn}

\begin{figure}[t]
    \centering
    \includegraphics[width = 0.5\linewidth]{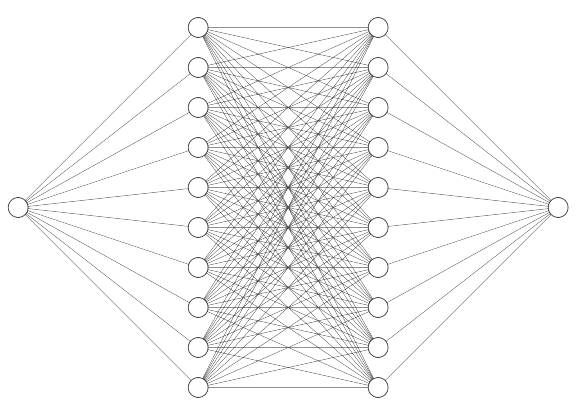}
    \caption{Image of a fully-connected neural network. }
    \label{fig:nn}
\end{figure}
In this section, we give a brief review of fully-connected neural networks (FCNNs)\cite{dnn} and PINNs. FCNN comprises the simplest type of neural networks in terms of structures. By the Universal Approximation Theorem \cite{hornik1989multilayer}, a FCNN with two hidden layers of infinite width is capable of approximating any function. In actuality, infinite width is unachievable. Networks with more hidden layers of large width are used as surrogate functions for approximations. 

A FCNN consists of an input layer, hidden layers, and an output layer. Within each layer, information passed in is processed by a composition of linear transformation with a weight and bias, followed by a nonlinear transformation $\sigma(\cdot)$. Suppose we have a FCNN with $L$ hidden layers, with input dimension $d_0$ and output dimension $d_{L+1}$. For each each layer, the input is the output of the previous layer and the output is of dimension $d_l$, for $l = 1, 2,\dots, L$. Let the collection of $N$ inputs $X_{\text{in}}$ be of dimension $d_0$, i.e. $X_{\text{in}}\in \R^{N\times d_0}$. 
% Suppose there is also a collection of $N$ inputs $X_{\text{in}}$ of dimension $d_0$, i.e. $X_{\text{in}}\in \R^{d_0\times N}$, and a nonlinear activation function $\sigma(\cdot)$.
A fully-connected network can be written in the following form:
\begin{equation}
    \begin{split}
        Y_1 &= \sigma(X_{\text{in}}W_1+b_1),\\
        Y_{l+1} &= \sigma(Y_lW_l+b_l), \text{ for }l = 1, 2, \dots, L-1,\\
        Y_{\text{out}} &= Y_LW_{L+1}+b_{L+1},
    \end{split}
    \label{eq:recurrent_net}
\end{equation}
where $Y_1$ is the output of the input layer, $Y_{\text{out}}$ is the output of the neural network. 
% \begin{equation}
% % \begin{split}
% X_{\text{output}} \in \R^{L+1} = W_{L+1}\sigma(W_L \sigma(W_{L-1}\cdots \sigma(xW_1_{\text{in}}+b_1)+\cdots+b_{L-1})+b_{L})+b_{L+1}.
% % Y_0 &= X_{\text{in}}\\
% % Y_{j+1} &= \sigma(W_jY_j+b_j), \quad j =  1, \dots, N\\
% % Y_{N+1} &= W_{N+1}Y_{N}+b_{N+1}
% % \end{split}
% \label{eq:nn}
% \end{equation}

Here, $W_{l}\in \R^{d_{l}\times d_{l+1}}$ are weight matrices and $b_l\in \R^{1\times d_{l+1}}$ are bias vectors in the $l$-th hidden layer, for $l = 1,\dots,L$.  Let $\btheta$ denote the collection of weights $W_{l}$ and biases $b_l$, for $l = 1,\dots, L+1$. 

Note that, from the point of view of linear algebra, the addition of bias in $l$-th hidden layer can be viewed as adding $\mathbbm{1}_{\text{in}}b_{l}$, where $\mathbbm{1}_{\text{in}}\in \R^{N\times 1}$ is a vector of all ones such that $\mathbbm{1}_{\text{in}}b_{l}\in \R^{N\times d_{l+1}}$ and the size $N$ is dependent on the input size $X_{\text{in}}$. In other words, the network can be written as
% \begin{align*}
% % \begin{split}
% X_{\text{output}} \in \R^{L+1} = W_{L+1}\sigma(W_L \sigma(W_{L-1}\cdots \sigma(xW_1_{\text{in}}+b_1\mathbbm{1}_{\text{in}})+\cdots+b_{L-1}\mathbbm{1}_{\text{in}})+b_{L}\mathbbm{1}_{\text{in}})+b_{L+1}\mathbbm{1}_{\text{in}}.
% % Y_0 &= X_{\text{in}}\\
% % Y_{j+1} &= \sigma(W_jY_j+b_j), \quad j =  1, \dots, N\\
% % Y_{N+1} &= W_{N+1}Y_{N}+b_{N+1}
% % \end{split}
% \end{align*}

\begin{equation}
    \begin{split}
        Y_1 &= \sigma(X_{\text{in}}W_1+\mathbbm{1}_{\text{in}}b_1),\\
        Y_{l+1} &= \sigma(Y_lW_l+\mathbbm{1}_{\text{in}}b_l), \text{ for }l = 1, 2, \dots, L-1,\\
        Y_{\text{out}} &= Y_LW_{L+1}+\mathbbm{1}_{\text{in}}b_{L+1},
    \end{split}
    \label{eq:recurrent_net_ll}
\end{equation}

To align notation with the following sections, as the network serves as a corrector, let $c_{\btheta}(x)$ denote such a network with inputs $X_{\text{in}}$. The physical significance of $x$ is discussed later.  In practice, the weight matrix $W_{l}$ in the $l$-th hidden layer is commonly initialized following the Xavier initialization \cite{glorot2010understanding}, where entries are independent and identically distributed (i.i.d.) drawn from a normal distribution $\mathcal N(0, \frac{2}{d_{l}+d_{l+1}})$. Bias vectors are commonly initialized to be zero. We will follow this convention in our experiments. 
% \subsection{Physics-Informed Neural Networks (PINNs)}
% Jagtap and Karniadakis proposed an adaptive activation function that accelerates convergence of PINN in \cite{JAGTAP2020109136}. For each of the layer in the neural network, instead of just learning the weights and bias in the operation $\sigma(Wx+b)$, the network learns an additional parameter $a$ in $\sigma(a(Wx+b))$. Jagtap claims that learning this addtional parameter assists the network in approximating the higher frequency components. We will test the effects of adaptive activation functions in our test problem. 

To introduce physics-informed neural networks (PINNs), we first need to formally specify a partial differential equation (PDE) that we are interested in solving. We are interested in solving a steady-state PDE
\begin{equation}
\begin{split}
\mathcal{N}[u](x) &= f(x), x\in \Omega,\\
u(x) &= g(x), x\in \partial\Omega,\\
\end{split}
\end{equation}
where $\mathcal N$ is a possibly nonlinear differential operator and $u(x):\bar{\Omega}\rightarrow \R$ is the solution with $x\in\R^{d_0}$. 
% Given an estimate $\hat{u}$ for the solution, the residual function is
% \begin{equation}
% r(x) = \mathcal{N}[\hat{u}](x)-f(x)
% \label{eq:residual_general}
% \end{equation}

Following \cite{RAISSI2019686}, we will explore the use of PINNs to compute corrections $c_{\theta}(x)$ to a given initial approximation $\hat{u}$, so that the surrogate has the form $u_s(x) = \hat{u}+c_{\theta}(x)$. The input is a set of spatial coordinates within the domain $\Omega$, and the output is $c_{\theta}(x)$. Different choices of $\hat{u}$ are discussed in Sections \ref{sec:correcting_solution} and \ref{sec:correcting_reduced}. The residual of the corrected solution $u_s(x)$ can be written as
\begin{equation}
r_{\btheta}(x) = \mathcal{N}[u_s(x)]-f(x)
\label{eq:theta_residual}
\end{equation}
% PINNs were first proposed in \cite{RAISSI2019686} as neural networks that utilize the governing equation, i.e., the underlying PDE, in the loss function. The idea of PINNs is to use the residual function $r$ to ensure the neural network follows the governing PDE. Let residual of a generic PDE be:

The network residual $r_{\btheta}$ can be derived through automatic differentiation \cite{van2018automatic}. The appropriate parameters $\btheta$ are obtained by minimizing a loss function:
\begin{equation}
\begin{split}
\mathcal L(\btheta) = \mathcal L_{u}(\btheta)+ \mathcal L_{r}(\btheta)\\
\mathcal L_{u}(\btheta) = \frac{1}{N_{u}}\sum_{i = 1}^{N_{u}}(u_s(x_{u}^{(i)})-u_{bc}^{(i)})))^2 \\
\mathcal L_{r}(\btheta) = \frac{1}{N_{r}}\sum_{i = 1}^{N_{r}}(r_{\btheta}(x_{r}^{(i)}))^2 
\end{split}
\label{eq:loss}
\end{equation}
where $\{(x_{u}^{(i)}),u_{bc}^{(i)}\}_{i=1}^{N_u}$ represent the values of training samples on the domain boundary, and $\{(x_{r}^{(i)})\}_{i=1}^{N_r}$ are points sampled within $\Omega$ for minimizing the residual. Here $N_u$ and $N_r$ are sample sizes for boundary conditions and residual points, respectively. 
Such innocent-looking mean squared errors lead to difficult training processes, as noted in \cite{wang2022and}. Understanding the training processes may lead to better designs of network architecture and loss functions for the desired governing equation. 

% \subsection{Input Transformation}
It is a common practice to normalize the inputs to a network as scaling of the inputs plays an important role in easing the backpropogation process \cite{normalizeinput}. This procedure can improve the training and testing results of the network remarkably, which lead us to wonder about the effects of input transformation. When the spatial coordinates in the training samples and the corresponding targets are fixed and only trainable weights and biases are being optimized, transforming the spatial coordinates leads to different loss functions, which could yield an easier training process and lead to more accurate approximations. The explicit form of loss functions will be explored in later sections. In this paper, we propose an input transformation scheme consisting of, for spatial coordinates $x$, a linear transformation $T(x)$ on the input coordinates $x$, 
\begin{equation}
T(x) = a(x+b),
\label{eq:transformation}
\end{equation}
so that input spatial coordinates $x$ in the training samples is replaced by $T(x)$.

% Let the new network structure be
% \begin{equation}
% {c}_{\btheta} = W_{L+1}\sigma(W_L \sigma(W_{L-1}\cdots \sigma(W_1(a(x+b))+b_1)+\cdots+b_{L-1})+b_{L})+b_{L+1}.
% \end{equation}
% where $a,b\in \R$.

\section{The Convection-Diffusion Equation}
\label{sec:cd}
\begin{figure}[h]
\centering
% \begin{subfigure}{0.48\textwidth}
% \centering
\includegraphics[width = 0.5\linewidth]{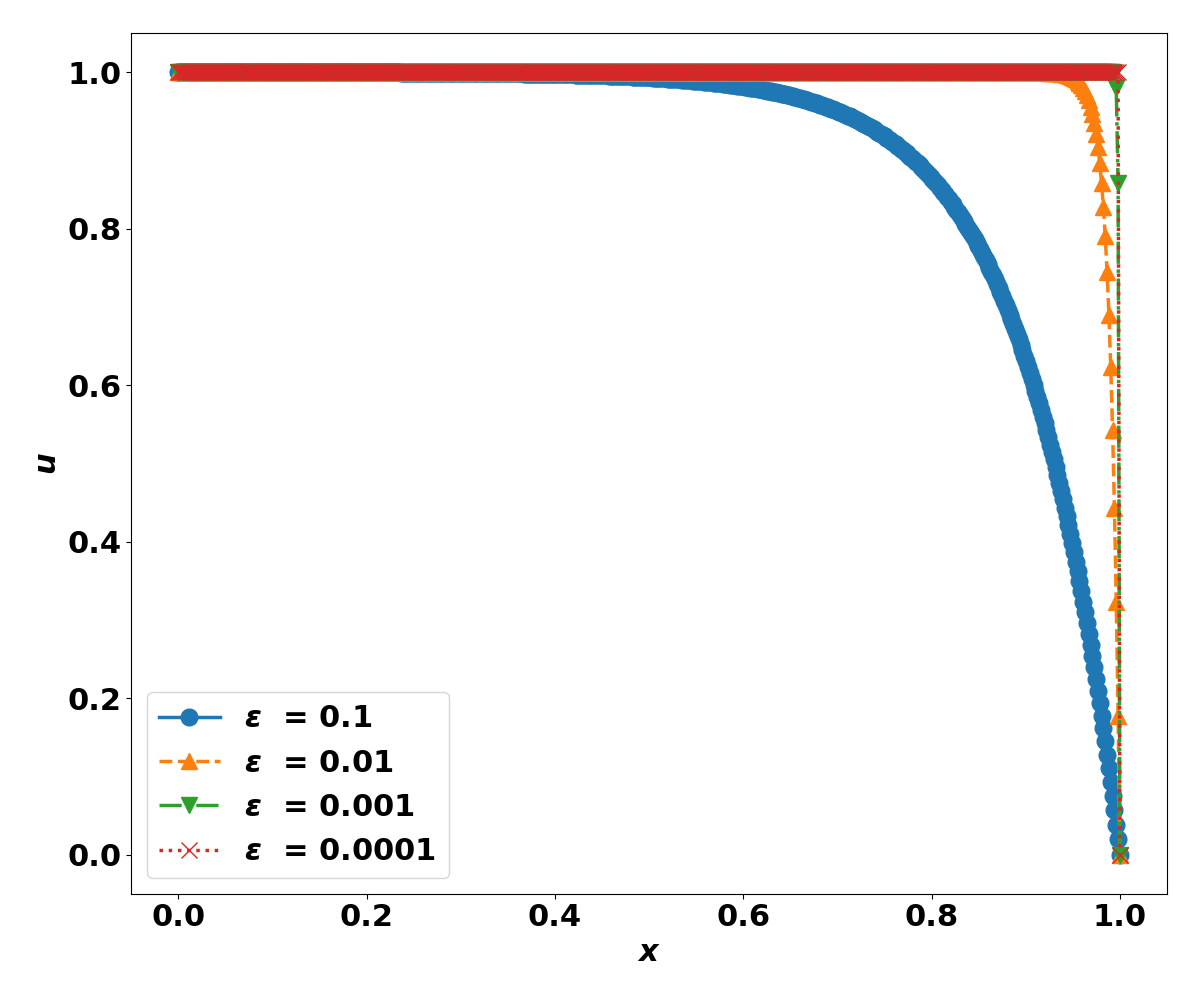}
% \caption{}
% \label{fig:1d_truesol}
% \end{subfigure}
\caption{Exact solutions of Equation \eqref{eq:1dcd2} for $\epsilon\in \{0.1,0.01,0.001,0.0001\}$.}
\label{fig:1dcd_truesols}
\end{figure}

In this section, we briefly describe the equations we will use as test cases. The steady-state convection-diffusion equation is given by
\begin{equation}
-\epsilon \nabla^2 u+\vec{b}\cdot \nabla u = f, \quad x\in \Omega \subseteq \R^d
\label{eq:cd}
\end{equation}
where $\epsilon>0$, $\nabla^2 u$ corresponds to diffusion, $\nabla u$ represents convection, $\vec{b}$ denotes a velocity and $f$ is the source term. The one-dimensional version is
\begin{equation}
\begin{split}
-\epsilon u''+b(x)u'&=f(x), \quad \text{for }x\in(0,1), 
% &u(0)=u(1)=0.
\end{split}
\end{equation}
Boundary conditions are specified at $x=0$ and $x=1$. We will use Dirichlet boundary conditions in the test cases. A simple version of \eqref{eq:cd}, 
\begin{equation}
\begin{split}
-\epsilon u^{''}+u^{'} &= 0 \quad \text{for }x\in (0,1),\\
u(0) &= 1-e^{-1/\epsilon},\quad u(1) = 0,
\end{split}
\label{eq:1dcd2}
\end{equation}
has the analytic solution 
\begin{equation}
u_{ex}(x,\epsilon) = 1-e^{(x-1)/\epsilon}.
\label{eq:1dcd2ex}
\end{equation} 
Figure \ref{fig:1dcd_truesols} shows these solutions for different values of $\epsilon$, where boundary layers form near $x=1$, and those layers become steeper as $\epsilon$ approaches $0$.

As is well known, these boundary layers cause the problem to be difficult to solve by standard methods. We will show below that they also affect the performance of PINNs, for example, leading in some cases to surrogate approximations to a \textit{different} convection-diffusion equation solution. We will investigate the cause of these tendencies in Sections \ref{sec:correcting_solution} and \ref{sec:correcting_reduced}. 

% The one-dimensional convection-diffusion equations form a good problem set for understanding the behavior of PINNs. On one hand, we could analytically solve for the exact solutions to these one-dimensional problems. The exact solution provide us a direct view on the accuracy of PINN-corrected solutions. Due to its steep boundary layer, the convection-diffusion equations serve as great test cases to investigate the behavior of PINNs under such numerically difficult cases. Later we observed that PINNs produced approximations that deviate from the exact solution to the governing equation. But the produced approximations are solutions to another convection-diffusion equation with the signs of the diffusion term inverted. The found steep boundary layer appeared in an incorrect place. We investigate the cause of such tendencies in later sections, through the help of input transformations. 

% We recognize that the equations used in real-world application are far more complicated than the problems we are about to discuss. However, as convection-diffusion equation in one-dimensional still presents itself as a tough problem for PINNs \cite{LP}, we will start with the one-dimensional case. 

% % \section{Application of PINNs as a Correction to Discrete Solutions}

\section{Correcting Oscillatory Discrete Solutions}
\label{sec:correcting_solution}
\begin{figure}[h]
\begin{subfigure}{0.48\textwidth}
\centering
\includegraphics[width = \linewidth]{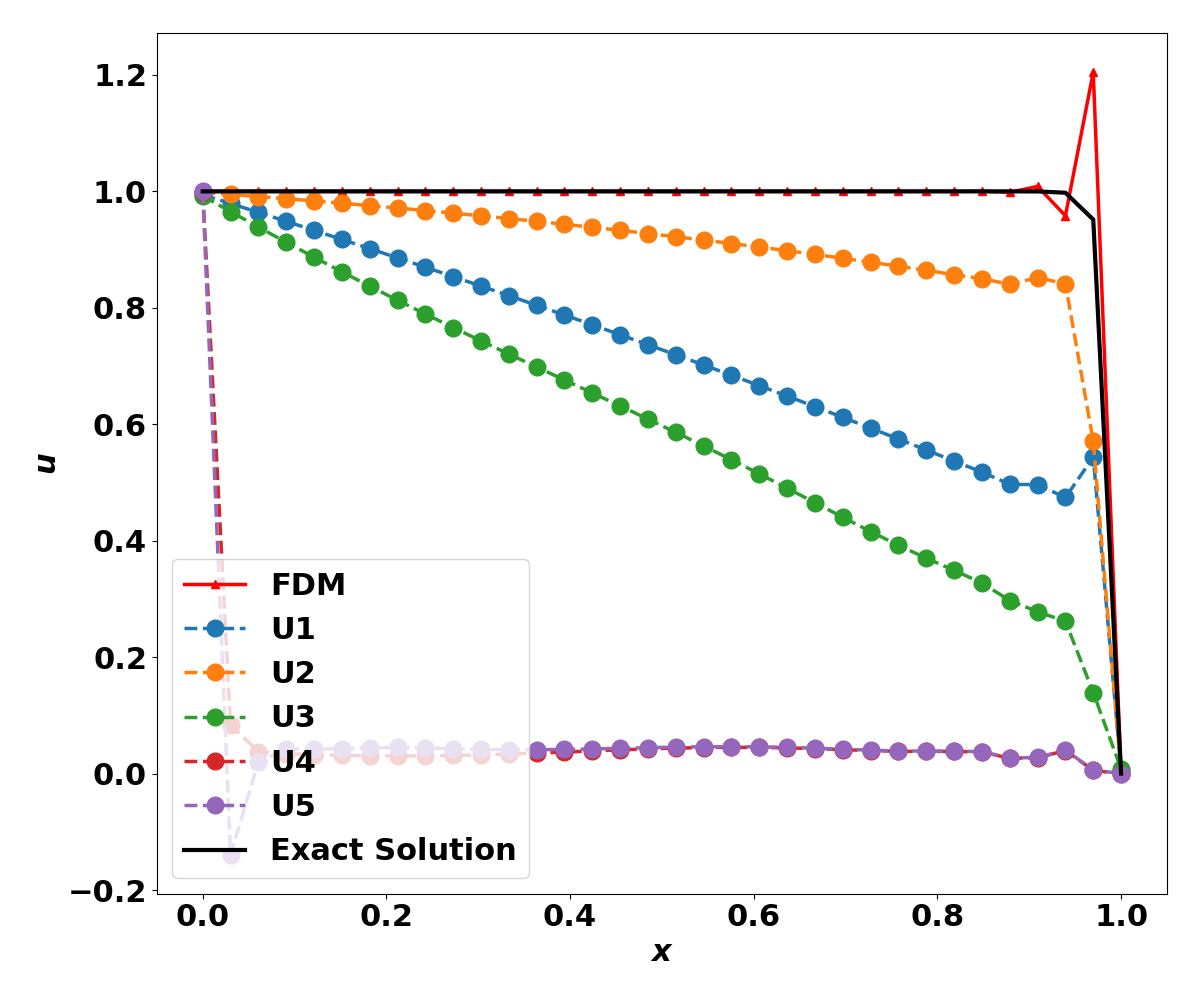}
\caption{}
\label{fig:1dcdd1e-2}
\end{subfigure}
\begin{subfigure}{0.48\textwidth}
\centering
\includegraphics[width = \linewidth]{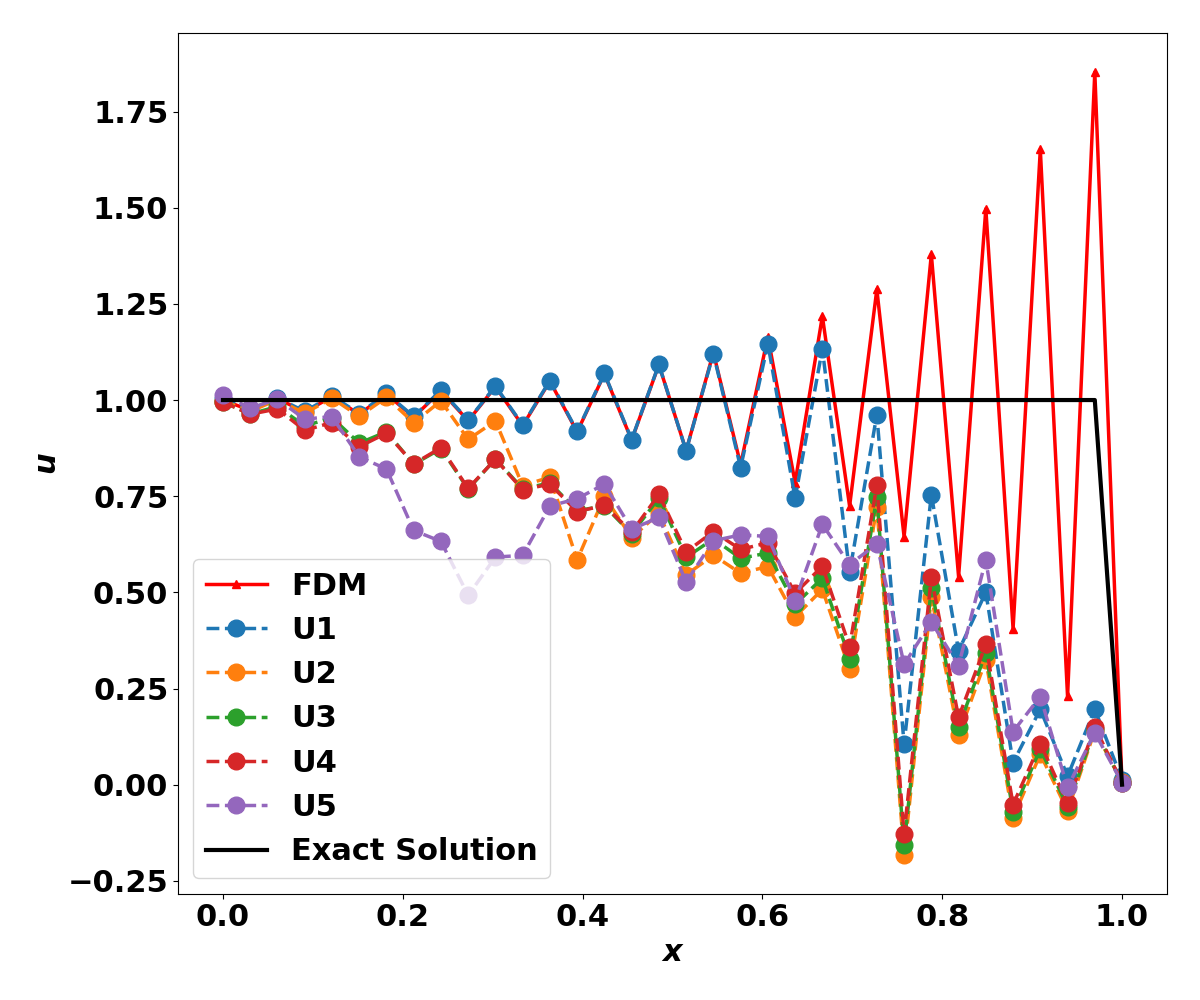}
\caption{}
\label{fig:1dcdd1e-3}
\end{subfigure}
\caption{Corrections made by standard PINN on FDM oscillatory solutions for: $\epsilon=0.01$ in Figure \ref{fig:1dcdd1e-2} and $\epsilon=0.001$ in Figure~\ref{fig:1dcdd1e-3}. }
\label{fig:1dcd_discrete}
\end{figure}
Since the introduction of PINNs, there have been attempts to use them to resolve boundary layers \cite{PATEL2022110754,Gasmi2021,fuks2020limitations}. Although the exact solution \eqref{eq:1dcd2ex} to the convection-diffusion problem is continuous, the solutions can be viewed as being close to discontinuous for small $\epsilon$, which poses a new issue on using FCNN as surrogate functions. As shown in \cite{PATEL2022110754, discontinuity}, neural networks struggle with approximating discontinuities. There is a necessity to investigate how FCNNs behave when PINN is used to solve problems with steep boundary layers. In this section, we examine the effectiveness of PINNs to correct oscillatory discrete solutions to the convection-diffusion equation obtained from FDMs. 

Standard PINN approaches treat FCNN as the surrogate function for approximating the solutions to PDE problems \cite{raissi2020hidden}. We took a different angle: instead of using FCNNs to directly approximate solutions, we use FCNNs as corrections to existing discrete solutions. Inspired by iterative numerical methods, where a new correction is applied to the current approximation at each iteration, we apply the same idea on using FCNNs as corrections. Starting at the initial iteration, we train an FCNN using a  method described below to correct the initial discrete solutions. If the corrected approximations do not achieve satisfactory accuracy, we train another FCNN to correct the corrected discrete solutions. Such iterations can be performed repeatedly until an adequate discrete approximation is obtained. Suppose one were to attempt to solve equation \eqref{eq:1dcd2} with standard discretizations, such as the central finite difference method (FDM). Let $\{x_j\}_{j=0}^{N}$ be an equi-spaced mesh on $[0,1]$, where $x_0 = 0$ and $x_{N} = 1$, with mesh $h$. The central FDM formulation is 
\begin{equation}
\Big(-\frac{\epsilon}{h^2}-\frac{1}{2h}\Big)U_{i-1}+\Big(2\frac{\epsilon}{h^2}\Big)U_i+\Big(\frac{1}{2h}-\frac{\epsilon}{h^2}\Big)U_{i+1} = 0,
\label{eq:central_fdm}
\end{equation}
where $U_i$ is the discrete approximation of $u(x)$ at $x = x_i$. Here $U(0) = 1-e^{-1/\epsilon}$ and $U(1) = 0$. On a mesh that is too coarse, the discrete solutions are oscillatory and non-physical \cite{esa}. Let $\hat{u}$ be the piecewise linear function that connects discrete values $\{U_i\}_{i=0}^N$ at $\{x_i\}_{i=0}^{N}$. When the mesh Peclet number $\mathcal P_h = \frac{h}{2\epsilon}$ of \eqref{eq:central_fdm} is larger than unity, the boundary layer cannot be resolved on the mesh and oscillations appear in the discrete solution. Examples are shown in Figure \ref{fig:1dcd_discrete}. We would like find a network $c_{\btheta}$ such that $u_s(x) = \hat{u}+c_{\btheta}$ more closely satisfies the governing equation, i.e., $u_s(x)\approx u(x)$. Substituting $u_s(x)$ into \eqref{eq:1dcd2}, the formulation becomes
\begin{equation}
\begin{split}
-\epsilon (c_{\btheta})^{''}+c_{\btheta}^{'} = \epsilon\hat{u}^{''}-\hat{u}^{'} &=-\hat{u}^{'},  \quad \text{for }x\in (0,1),\\
u_s(0) = (\hat{u}+c_{\btheta})(0) &= 1-e^{-1/\epsilon},\\
u_s(1) = (\hat{u}+c_{\btheta})(1) &= 0.
\end{split}
\end{equation}
Note that $\hat{u}^{''}=0$ between $\{x_i\}_{i=0}^N$ since $\hat{u}$ is assumed to be piecewise linear. For points in between $\{x_i\}_{i=0}^N$, $\hat{u}^{'}$ are the slopes between $\{U_i\}_{i=0}^N$. When we select sample points for computing the residual mean squared error in the loss function, we avoid $\{x_i\}_{i=0}^N$ to circumvent complications with discontinuities. 

The loss function for $c_{\btheta}$ has the form of equation \eqref{eq:loss} except the residual is written as
\begin{equation}
r_{\btheta}(x_{r}^{(i)}) = -\epsilon (c_{\btheta})^{''}+c_{\btheta}^{'} + \hat{u}^{'}
\label{eq:loss_cd}
\end{equation}
With given initial approximations $\{U^0_i\}_{i=0}^{N}$, we can train $c_{\btheta}^j$ for correction iteration $j = 1, 2, \dots$. Each subsequent $c_{\btheta}^{j+1}$ will be a correction for $\hat{u}+\sum_{k=1}^j c_{\btheta}^k$. These correction steps can be performed until $(\hat{u}+\sum_{k=1}^j c_{\btheta}^k)(x)$ is smooth and accurately approximates $u(x)$. 
% Let the discrete solution $\{U_i^j\}_{i=0}^{N}$ denote $(\hat{u}+\sum_{k=1}^j c_{\btheta}^k)(x)$ evaluated at $\{x_i\}_{i=0}^N$. 
This process can be written recursively as 
\begin{equation}
(U^j)(x_i) = (U^{j-1}+c_{\btheta}^{j})(x_i), \quad \text{for }j = 1,2, \dots, i = 0,\dots, N,
\label{eq:recursive}
\end{equation} 
\textcolor{black}{where $\{U_i^j\}_{i=0}^{N}$ denote $(\hat{u}+\sum_{k=1}^j c_{\btheta}^k)(x_i)$. }
\begin{figure}[t]
\begin{subfigure}{0.48\textwidth}
\centering
\includegraphics[width = \linewidth]{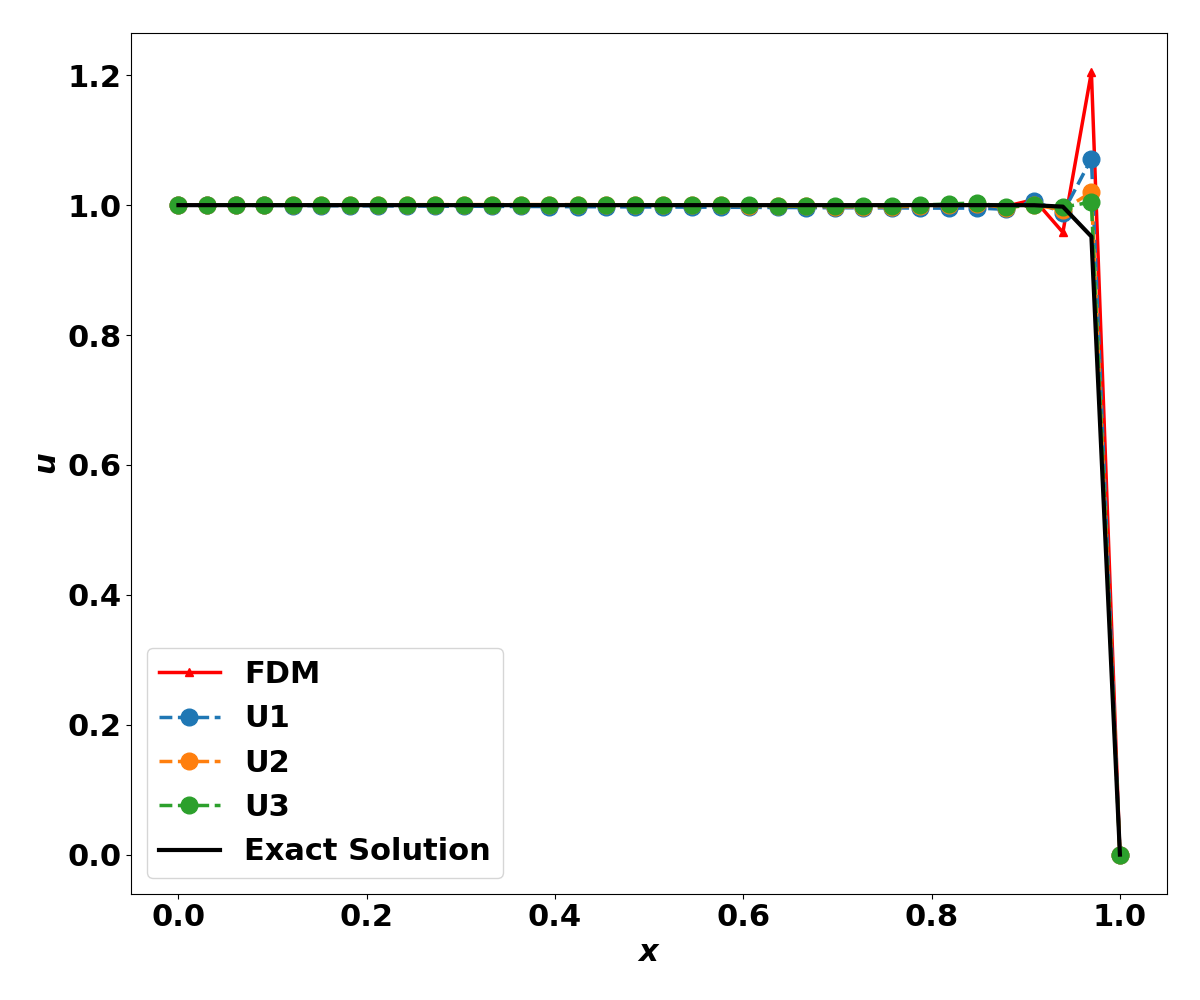}
\caption{}
\label{fig:1dcdd1e-2_trans}
\end{subfigure}
\begin{subfigure}{0.48\textwidth}
\centering
\includegraphics[width = \linewidth]{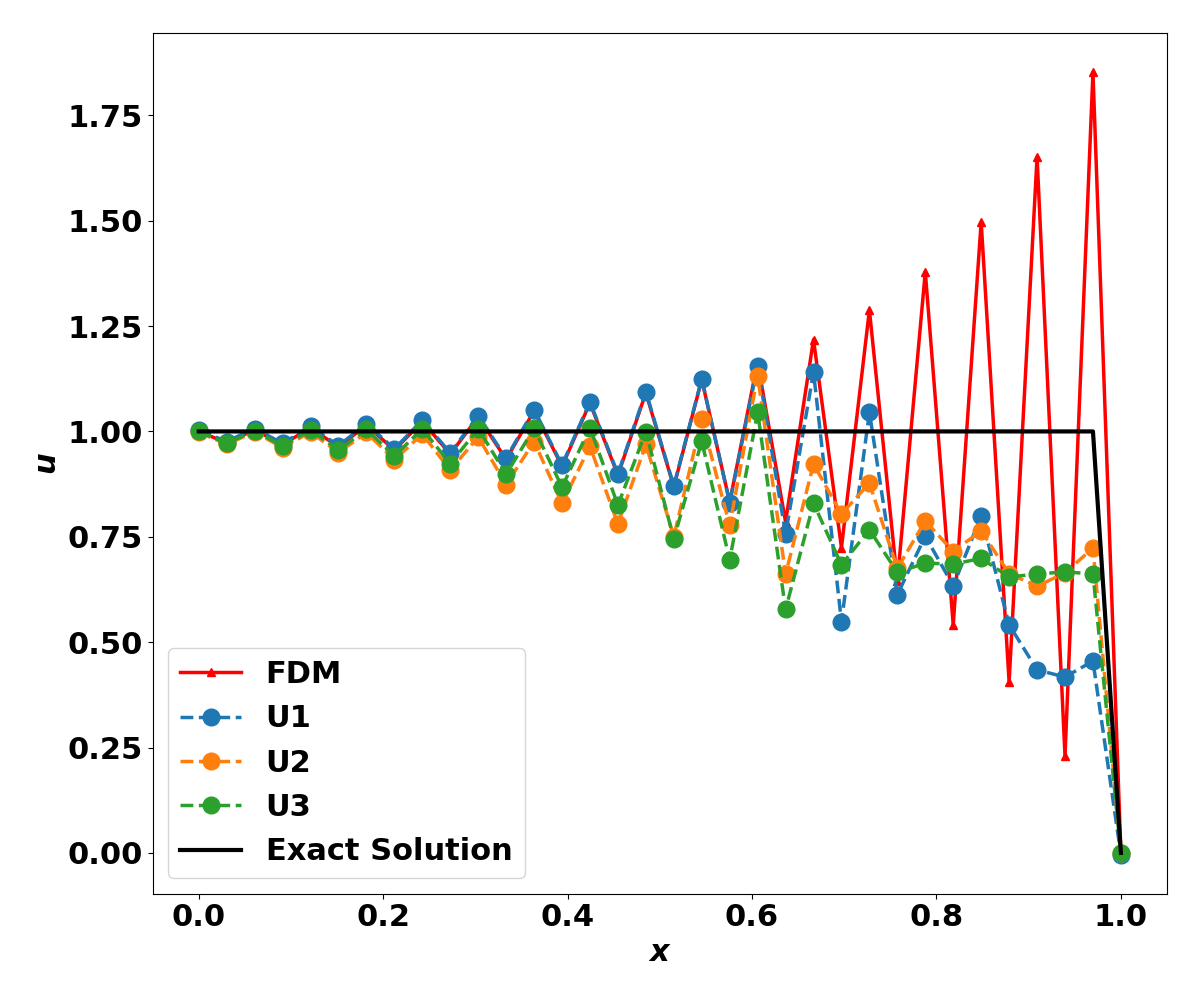}
\caption{}
\label{fig:1dcdd1e-3_trans}
\end{subfigure}
\caption{Corrections made by PINN with input transformation $T(x) =  100(x-1)$ on FDM oscillatory solutions for: $\epsilon=0.01$ in Figure \ref{fig:1dcdd1e-2_trans} and $\epsilon=0.001$ Figure~\ref{fig:1dcdd1e-3_trans}.}
\label{fig:1dcd_discrete_trans}
\end{figure}
\textcolor{black}{Representative examples of this process are shown} in Figure \ref{fig:1dcd_discrete} for $\epsilon = 0.01, 0.001$, on a uniform grid of size $N = 32$ over $x\in [0,1]$. The mesh Peclet numbers are larger than unity for both values of $\epsilon$, and the FDM solutions are oscillatory. An equi-spaced collection of points over $[0,1]$ with stepsize $1/128$ was used as the training residual sample collection. The uniform grid points used to generate FDM approximations \textcolor{black}{were not used as} sample points, \textcolor{black}{which rendered 96 points in the residual training samples.} We train individual FCNNs of 3 hidden layers with width 40 per layer using 10,000 epochs of the Adam optimization scheme for 5 correction iterations. For $\epsilon =0.001$, by the fifth correction, the approximations $U^5$ are \textcolor{black}{still far from the exact solution.  The corrections do not converge monotonically to an accurate surrogate function. Similarly, for $\epsilon=0.01$, approximations of $U^5$ present a flow opposite to the exact solution. } 

As discussed in Section \ref{sec:nn}, input transformations result in different loss functions. We explore the impact of transformations on this process, using the transformation $T(x) = 100(x-1)$ on inputs of the network. The resulting corrected approximation becomes $u_s(T(x)) = \hat{u}(x)+c_{\btheta}(T(x))$. Using the same network and training settings as in the previous tests, we found that $c_{\btheta}(T(x))$ performs better than $c_{\btheta}(x)$ for $\epsilon = 0.01$. This can be seen from Figure \ref{fig:1dcd_discrete_trans} where the corrected solution approaches the exact solution with each subsequent iteration for $\epsilon = 0.01$. By the third correction iteration, the corrected solution $U^3$ closely follows the exact solution. Unfortunately, this transformation is not effective for $\epsilon = 0.001$. The corrected solutions failed to capture the steep boundary layer. Instead of correcting to a steep boundary layer near $x=1$ around $y=1$, the corrected approximations develops an incorrect flat region near $y=0.75$. These results are representative of our experience for correcting oscillatory solutions, \textcolor{black}{and we found that using such transformations to correct oscillatory discrete solutions does not lead to a robust solution strategy.}
In the next section, we explore an alternative approach, where PINNs are used to correct a simpler solution, derived from the reduced equation. 
% Since the corrected approximations from oscillatory FDM solutions do not converge to the exact solution with more iterations of correction and the corrected approximations appear to deviate from the exact solutions, correcting discrete oscillations may not be a suitable task for standard PINNs. Oscillations might be too complicated for standard PINNs to correct.  

\section{Correcting Reduced Solutions}
\label{sec:correcting_reduced}
\begin{figure}[t]
\centering
\begin{subfigure}{0.48\linewidth}
\includegraphics[width = \linewidth]{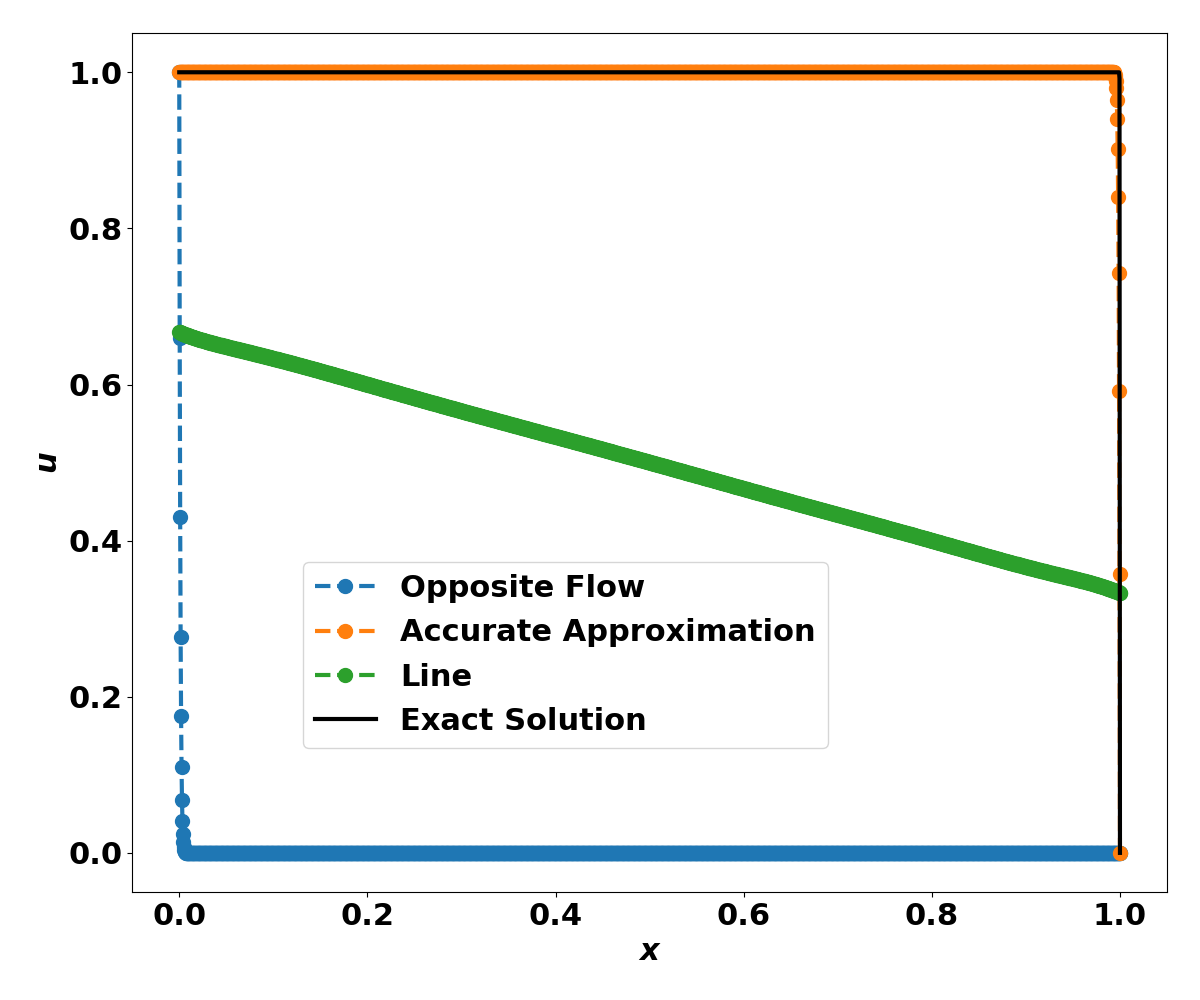}
\caption{}
% \caption{Different approximations generated using a standard PINNs under different random initialization seeds.}
\label{fig:seeds}
\end{subfigure}
\begin{subfigure}{0.48\linewidth}
\centering
\includegraphics[width = \linewidth]{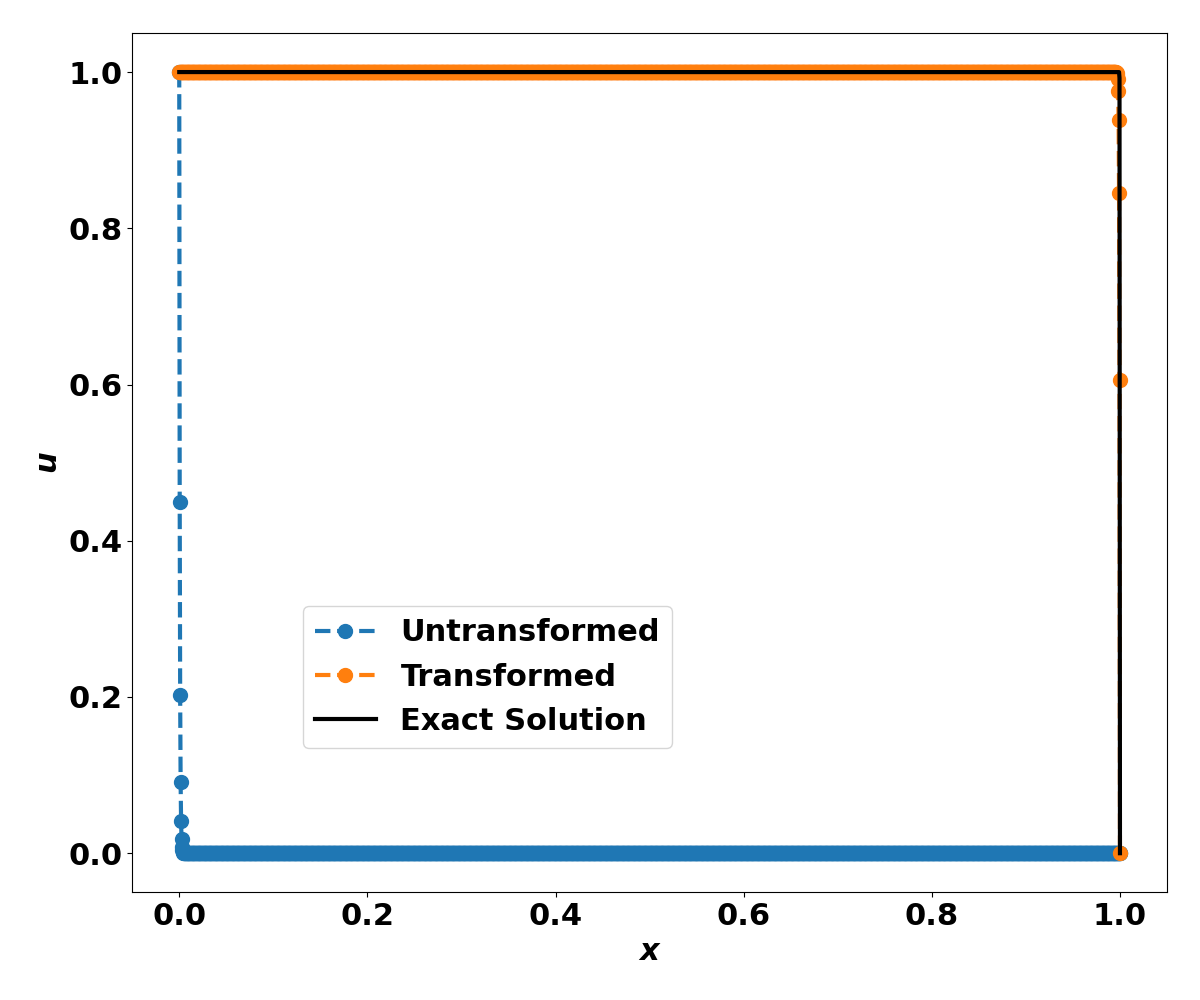}
\caption{}
\label{fig:ib}
\end{subfigure}
\caption{Different approximations generated using a standard PINNs under different random initialization seeds in \ref{fig:seeds}. Comparison of effects of transformation in \ref{fig:ib}. Here the input is transformed from $x\in [0,1]$ to $x-1\in [-1,0]$.}
\label{fig:reduced_solutions}
\end{figure}

%why correct reduced solution
\textcolor{black}{In this section, we consider an alternative approach for beginning the PINN procedure.} Instead of the test problem in \eqref{eq:1dcd2}, consider the following \textit{reduced problem}
\begin{equation}
\begin{split}
u_0' &= 0,\\
u_0(0) &= 1,
\end{split}
\label{eq:reduced}
\end{equation}
which corresponds to the hyperbolic problem arising from $\epsilon=0$. The reduced solution is $u_0(x) = 1$ for $x\in [0,1]$, and there is no boundary condition (or layer) at the outflow boundary $x=1$. In this section, we take a different point of view and apply PINNs as corrections to the reduced solutions of the convection-diffusion equations, i.e., the solution to the pure convection equation for $\epsilon = 0$. \textcolor{black}{We set $\hat{u}(x) = u_0(x)$. The governing equation that $c_{\btheta}$ needs to satisfy is 
\begin{equation}
    -\epsilon c_{\btheta}''+ c_{\btheta}' = 0
    \label{eq:reduced_governeqn}
\end{equation}
}
\subsection{Input Transformation: Shifting}

\textcolor{black}{We first explore experimentally} the effect of shifting inputs on PINN correcting reduced solutions of the convection-diffusion equation. We would like to train $c_{\btheta}$ such that $u_0+ c_{\btheta}$ approximates the solution to \eqref{eq:1dcd2}.\footnote{Different one-dimensional convection-diffusion equations may need similar correction of different reduced solutions. An example is 
\begin{equation}
\begin{split}
-\epsilon u^{''}+u^{'} &=  1,\\
u(0) = u(1) &= 0.
\end{split}
\label{eq:1dcd_comp}
\end{equation}
The reduced solution of \eqref{eq:1dcd_comp} is $x$. But the correction needed satisfies the residual function of $-\epsilon c_{\btheta}^{''}+c_{\btheta}^{'} =  0$ with boundary conditions $c_{\btheta}(0) = 0$ and $c_{\btheta}(1) = -1$. The correction needed for \eqref{eq:1dcd_comp} has similar characteristics to the one for \eqref{eq:1dcd2}. Both corrections contain a steep gradient layer near $x=1$. } 

In one initial attempt to approximate the correction to the reduced solution, we treated $c_{\btheta{}}$ using a standard PINN approach. For residual samples, we create an equi-spaced mesh over $x\in (0,1)$ with stepsize $h=1/128$, \textcolor{black}{as we did} for correcting FDM solutions. The boundary samples are $\{(0, -e^{-1/\epsilon}), (1,-1)\}$. The loss function is the standard PINN loss function in \eqref{eq:loss}. For a network of one hidden layer, width 20, we train the network with Adam optimization \cite{kingma2017adam} for 10000 epochs, followed by 10000 epochs of LBFGS \cite{liu1989limited}. \textcolor{black}{
Weights and biases of the networks were initialized following the Xavier initialization using different random seeds. Using 100 different random seeds, we found that only \textcolor{black}{seven} trained networks produced good approximations that closely followed the exact solution. \textcolor{black}{Eighty-six} of the trained networks produced approximations with a steep boundary layer located near $x=0$, corresponding to a flow in the direction opposite to that of the solution, as for the solution of $-\epsilon u''-u' = 0$ with the same boundary conditions as \eqref{eq:1dcd2}. The remaining \textcolor{black}{seven} of the trained networks produced lines with negative slopes. The three types of approximations are shown in Figure \ref{fig:seeds}. } \textcolor{black}{We also found that increasing the depth or width of the network does not improve the chance for PINNs obtaining an accurate approximation.} 
% As shown in Figure \ref{fig:seeds}, there are three types of approximations the trained network found. The first one is a line between the boundary conditions. The second is the solution of $-\epsilon u''-u' = 0$ with the same boundary conditions as \eqref{eq:1dcd2}. The last is the correct correction to the reduced solution of \eqref{eq:1dcd2}. 

This points to the questions of whether it is possible for a network to reliably obtain good approximations and the cause of such different results. Our intuition was that inaccurate surrogates can be attributed to the properties of the activation function used, which was $\tanh(\cdot)$. The hyperbolic tangent function plateaus as the independent variable approaches $\pm \infty$, and its first derivative $\sech^2(\cdot)$ attains the largest value at the origin. 
\textcolor{black}{To approximate the steep boundary layer, we wondered if shifting the input of the network so that the location of the boundary layer aligns with root of the inputs of $\tanh(\cdot)$ could aid the training. The aim is to shift the inputs such that $\tanh{(\cdot)}$ attains its largest first derivative at the location of the boundary layer.} For input of the network $x\in [0,1]$, the input of the activation function of the first hidden layer is $xW_1+b_1$. 
As $\tanh(\cdot)$ attains its steepest gradient near the origin, it seemed beneficial to align the root of $xW_1+b_1$ with the location of the boundary layer at the initialization of the network. As described in section \ref{sec:nn}, weights are initialized following a normal distribution of mean zero and biases are initialized to be zeros. With this convention, $xW_1+b_1 = xW_1$. The root of $xW_1+b_1$ at network initialization is expected to be $x=0$, which does not naturally align with the location of the boundary layer $x=1$. Shifting the input to the right by 1 shifts the root of $(x-1)W_1+b_1$ to $x=1$, i.e. the location of the boundary layer. We performed the same training procedure as previous tests after shifting the inputs and tested on the same 100 random seeds. \textcolor{black}{With the shifted inputs, only three approximations corresponded to the flow in the wrong direction. The rest of the ninety-seven approximations are close to the exact solution. Note that we perform input shifting with prior knowledge about where the boundary layer is. For a simple one-dimensional case, the location of the boundary layer always appears near the outflow boundary. Such prior information is not difficult to obtain. More complex cases where the locations of the boundary layers are not known a priori are not discussed here.  }

We denote networks trained with unshifted inputs as \textit{Untransformed} networks and the ones trained with shifted inputs as \textit{Transformed} networks. The opposite flow approximation produced by an \textit{Untransformed} network is shown in Figure \ref{fig:ib} along with the accurate approximation produced by a \textit{Transformed} network. 
% The approximation produced by the network trained with shifted inputs is shown as \textit{Transformed} in the same figure. 
% Finding the correct approximations does not depend on the random initialization seeds. 
% \textcolor{black}{Both \textit{Transformed} and \textit{Untransformed} networks were of one hidden layer with width 2, and trained using 10,000 Adam epochs and 10,000 LBFGS epochs with the same training samples. }
To explain the drastically different performance obtained by the two types of networks, we first took a look at the difference in trained weights and biases between these two cases. 
Explaining the behavior of large network tends to be challenging due to the large number of parameters. However, with a small network, understanding performance is more tractable. Thus, we began our exploration with \textit{Transformed} and \textit{Untransformed} networks of one hidden layer with width 2. 

% Explaining the behavior of networks has always been a challenging task since number of parameters involved increases exponentially as more hidden layers are used. Luckily, we can reduce the network structure to one hidden layer with width 2 and still solve the problem. We first examine the difference in weights and biases under this small network structure. 
% Now we consider the transformed and untransformed networks. For the untransformed network, we have:
% $$\hat u_{\btheta}(x) = W_2\tanh(xW_1+b_1)+b_2$$
% For the transformed inputs case, the networks is
% $$\tilde u_{\btheta}(x) = W_2\tanh((x-1)W_1+b_1)+b_2$$
\begin{figure}[tp]
\centering
\begin{subfigure}{0.48\textwidth}
\centering
\includegraphics[width = \linewidth]{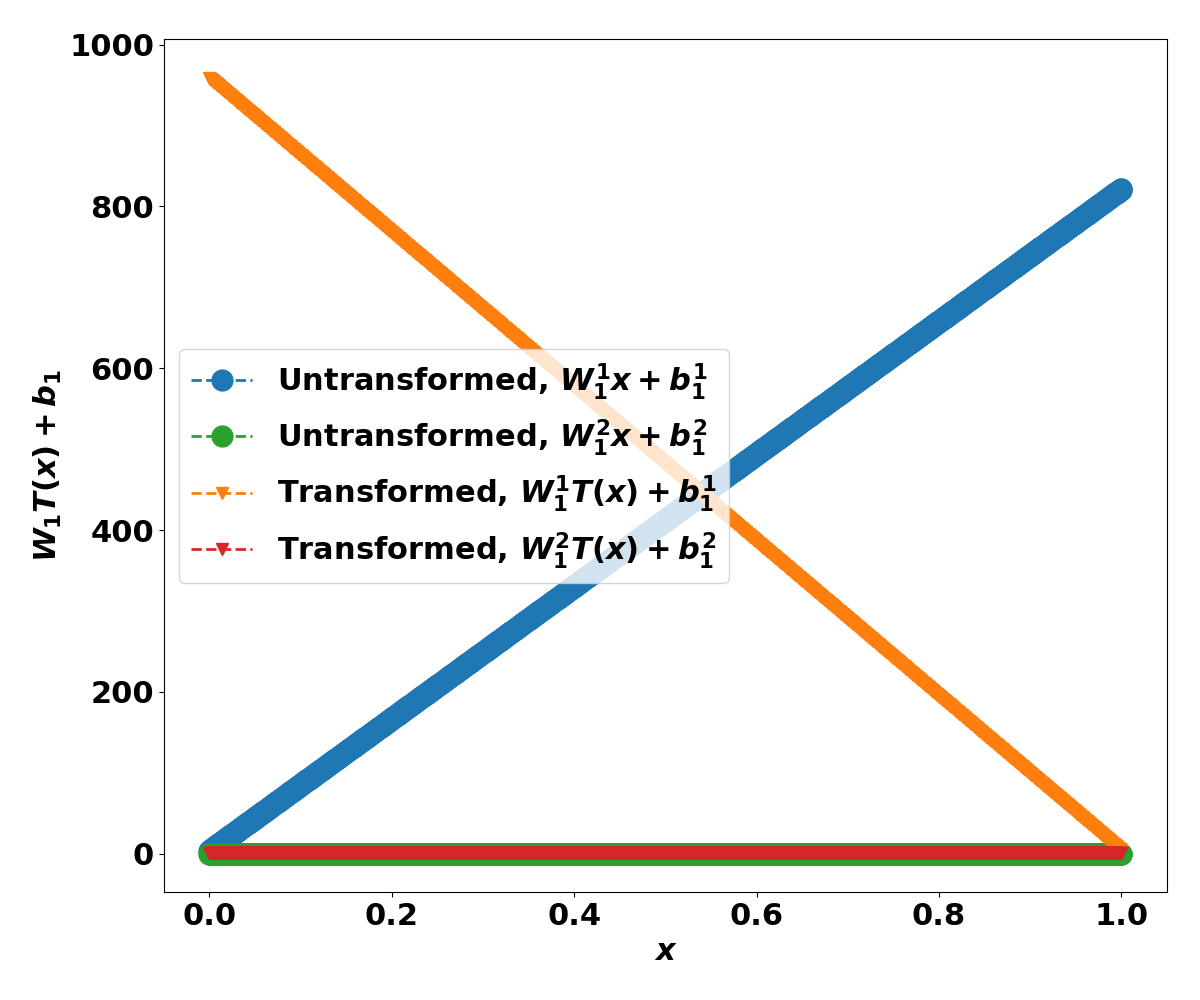}
\caption{Components of $T(x)W_1+b_1$}
\label{fig:afteraffine}
\end{subfigure}
\begin{subfigure}{0.48\textwidth}
\centering
\includegraphics[width = \linewidth]{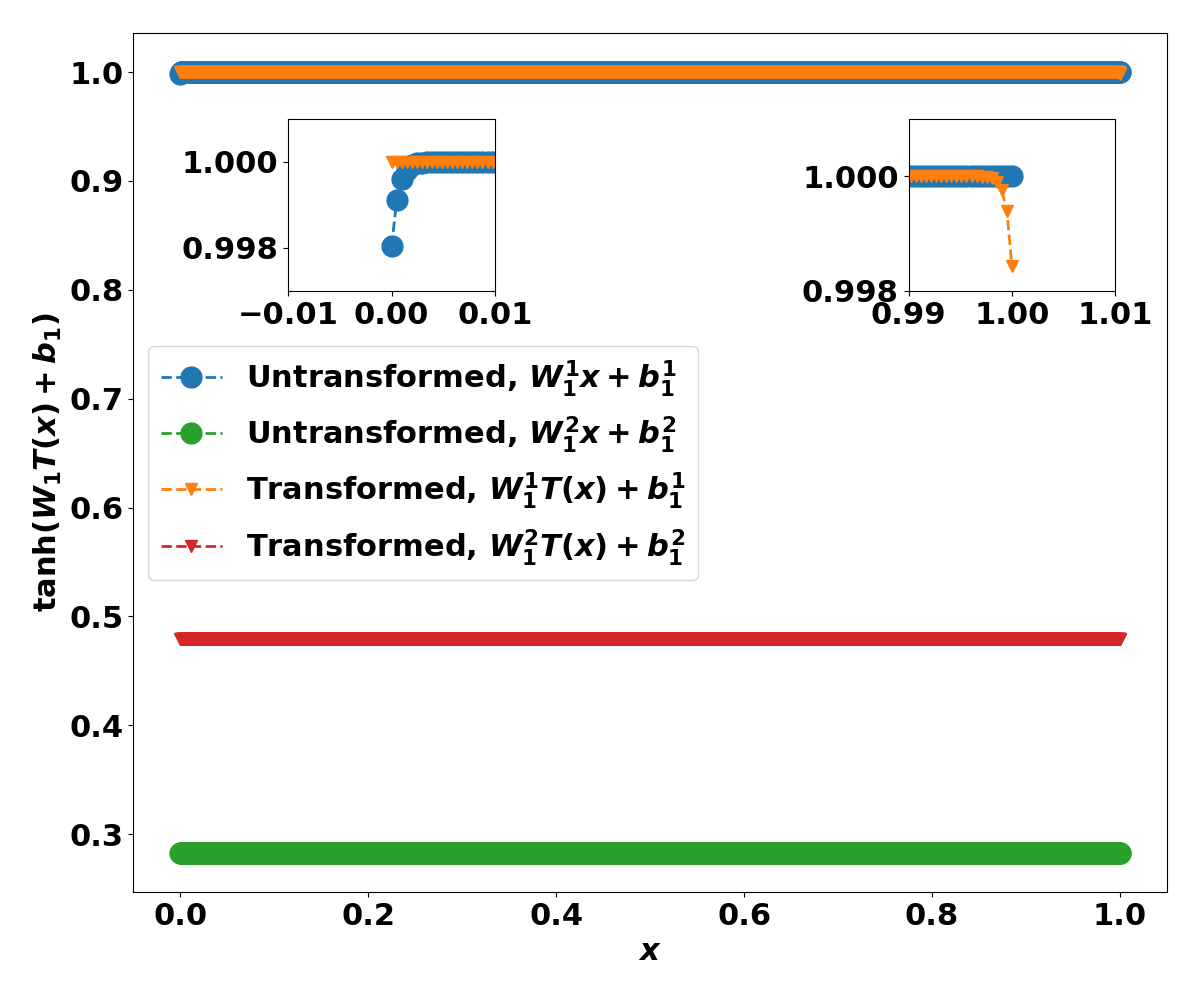}
\caption{Components of $\tanh(T(x)W_1+b_1)$}
\label{fig:afteractivation}
\end{subfigure}
% \centering
% \begin{subfigure}{0.96\textwidth}
% \centering
% \includegraphics[width = \linewidth]{Figures/after_activation_zoom.png}
% \caption{Zoomed in near 1 of Figure \ref{fig:afteractivation}. }
% \label{fig:zoomedin}
% \end{subfigure}
\caption{For a network with one hidden layer of width 2, in the first hidden layer, the output after the affine transformation is $xW_1+b_1$ (or $(x-1)W_1+b_1$), a matrix of two columns. We plot the two columns of the output of the affine transformation in \ref{fig:afteraffine}. We plot the two columns of the output of $\tanh(xW_1+b_1)$ (or $\tanh((x-1)W_1+b_1)$) in \ref{fig:afteractivation}. Two inset figures show the zoomed-in look near $y=1$ at $x=0$ and $x=1$. 
% A zoomed-in look of Figure \ref{fig:afteractivation} near $x=1$ is shown in \ref{fig:zoomedin}. 
}. 
\label{fig:firstlayer}
\end{figure}

Inputs of the network $X_{\text{in}}$ come in a column vector of size $N\times 1$, \textcolor{black}{where $N$ is the number of input samples}. For a network of one hidden layer of width $d$, we have $W_1\in \R^{1\times d}$, $b_1\in \R^{1\times d}$, $W_2\in \R^{d\times 1}$, $b_2\in \R$. The network can be written as 
\begin{equation}
c_{\btheta}(X_{\text{in}}) = \tanh(X_{\text{in}}W_1+\mathbbm{1}_{\text{in}}b_1)W_2+b_2.
\label{eq:onelayernet}
\end{equation}
The output $c_{\btheta}(X_{\text{in}})$ is of size $\R^{N\times 1}$. 
% We apply training as previously stated for untransformed and transformed inputs.  
Figure \ref{fig:firstlayer} shows the output of the first layer of the trained networks. The correction obtained from the transformed input with transformation $T(X_{\text{in}}) = X_{\text{in}}-\mathbbm{1}_{\text{in}}$ is denoted $c_{\btheta}(T(X_{\text{in}})) = c_{\btheta}(X_{\text{in}}-\mathbbm{1}_{\text{in}}) =\tanh((X_{\text{in}}-\mathbbm{1}_{\text{in}})W_1+\mathbbm{1}_{\text{in}}b_1) W_2+\mathbbm{1}_{\text{in}}b_2$. 
% We first thought these two cases are initialized to have different expected loss values. However, the expected values of the loss of the two cases are the same. Then we start to consider how the weights and bias look after training. 

After training, the weight $W_1$ consists of one large-magnitude entry and one near-zero entry for both the transformed and untransformed cases. \textcolor{black}{For example, for the transformed case,  $W_1^1x+b_1^1$ is a line with a negative slope that varies from 0 to around 900 in a span of $x\in [0,1]$, while $W_1^2x+b_1^2$ is a line with a slope near zero that stays close to the $x$ axis. }
The large-magnitude entry in $W_1$ leads to a line with a steep slope after the affine transformation, as seen in Figure \ref{fig:afteraffine}. This figure shows components of the untransformed quantities $xW_1+b_1$ and transformed ones $T(x)W_1+b_1$. The steep lines cause the activation function $\tanh(\cdot)$ to plateau near 1 (or -1) for a majority of the domain $x\in [0,1]$. At either one of the boundaries, the line $xW_1+b_1$ (or $T(x)W_1+b_1$) approaches zero. The boundary closer to the root of the line manifests itself as a steep gradient after the activation function $\tanh(\cdot)$. The two inset figures in Figure \ref{fig:afteractivation} show a zoomed-in look of $\tanh(xW_1+b_1)$ and $\tanh(T(x)W_1+b_1)$ near the largest values on the $y$ axis at $x=0$ and $x=1$. 
\textcolor{black}{The gradients near the two boundaries are stretched into steep gradients shown in Figure \ref{fig:ib} after passing through the affine transformation in the output layer. Note that the affine transformation in the output layer does not change the character of $\tanh(T(x)W_1+b_1)$. Steep boundary layers cannot form if gradients are not seen in $\tanh(T(x)W_1+b_1)$. }
% It can be seen that the transformed case $\tanh(T(x)W_1+b_1)$ forms a steep gradient near $x=1$. 
% \footnote{$xW_1+b_1$ for the untransformed case and $T(x)W_1+b_1$ for the transformed case}

As noted previously, during training, all samples are fixed and only the weights and biases are being trained. The difference in inputs of samples leads to different loss functions, which in turn lead to different surrogate functions. To illustrate this, consider the analytic form of the loss function when the network size is fixed to a FCNN of one hidden layer. The goal of examining the analytic form is to find out the preferable PINN setups for the convection-diffusion equation. 
\textcolor{black}{Note that the residual of $c_{\btheta}$ is 
\begin{equation}
     -\epsilon c_{\btheta}''+ c_{\btheta}' = r.
    \label{eq:reduced_residual}
\end{equation}
When the residual $r=0$, 
% when the governing equation \eqref{eq:1dcd2} is satisfied, the residual function is zero, i.e.  
% $$-\epsilon u^{''}+u' = 0. $$
% Consequently, 
$c_{\btheta}$ satisfies }
\begin{equation}
\frac{c_{\btheta}^{''}}{c_{\btheta}'}=\frac{1}{\epsilon}.
\label{eq:ration}
\end{equation}

\textcolor{black}{
% A perfectly trained network $c_{\btheta}$ understand how the network $c_{\btheta}$ behaves, 
Consider the simplest network with one hidden layer and width one, i.e. $d=1$ in \eqref{eq:onelayernet}. 
% $${c}_{\btheta}(x) = W_2\tanh(x W_1 + b_1)+b_2$$
The derivatives of $c_{\btheta}$ with respect to $x$ are straightforward in this case, for $x\in [0,1]$, 
\begin{align*}
{c}_{\btheta}' &= W_2W_1\sech^2\left(x W_1 + b_1\right),\\
{c}_{\btheta}^{''} &= -2W_2(W_1)^2\tanh(x W_1 + b_1)\sech^2(x W_1 + b_1).
\end{align*}
Therefore, if $c_{\btheta}'\neq 0$,
\begin{align*}
\frac{{c}_{\btheta}^{''}}{{c}_{\btheta}'} = -2 W_1\tanh(x W_1 + b_1).
\end{align*}
% For the network $c_{\btheta}$ to minimize the residual function, $c_{\btheta}$ should satisfy $\frac{{c}_{\btheta}^{''}}{{c}_{\btheta}'} = \frac{1}{\epsilon}$. 
Thus, by \eqref{eq:ration},  
$$-2 W_1\tanh(x W_1 + b_1) = \frac{1}{\epsilon}.$$
Since $-1<\tanh(x)<1$ for $x\in \R$, $ \mathopen|-2W_1\tanh(x W_1 + b_1)\mathclose|<2\mathopen|W_1\mathclose|$. Thus, by triangle inequality, 
\begin{equation}
2\mathopen|W_1\mathclose|> \frac{1}{\epsilon},
\label{eq:ineq1}
\end{equation}
as $\epsilon>0$. We can infer from this inequality that $W_1$ should be very large in magnitude when $\frac{1}{\epsilon}$ is large, which coincides with our observation in Figure \ref{fig:firstlayer}, \textcolor{black}{where one entry in $W_1$ is large in magnitude}. For $\epsilon = 10^{-4}$, $\mathopen|W_1^1\mathclose|$ is around $10^3$ for both the transformed and untransformed cases.  \\
Now that we understand how the simplest network behaves, we can increase the width while fixing the number of hidden layers. }
Consider the network with one hidden layer of width larger than one ($d_1>1$), for $x\in [0,1]$, 
$${c}_{\btheta}(x) = \sum_{i = 1}^{d_1} W_2^i\tanh(x W_1^i + b_1^i)+b_2,$$
where the subscripts indicate the layer index and the superscripts indicate the entry index. The entry-wise explicit form is
\begin{align*}
{c}_{\btheta}' &=\sum_{i = 1}^{d_1} W_2^i \sech^2(x W_1^i + b_1^i)(W_1^i).\\
% {c}_{\btheta}^{''} &= -2\sum_{i = 1}^{d_1} W_2^i\tanh(x W_1^i + b_1^i)\sech^2(x W_1^i + b_1^i){(W_1^i)}^2.
\end{align*}
% \textcolor{black}{
For $(c_{\btheta}+u_0)(x)$ that satisfies \eqref{eq:1dcd2}, $c_{\btheta}$ should follow
\begin{equation}
    \begin{split}
        {c}_{\btheta}'(1) &= -\frac{1}{\epsilon}.\\
        % {c}_{\btheta}''(x=1) &= -\frac{1}{\epsilon^2}.
    \end{split}
    \label{eq:deri_c}
\end{equation}
Therefore, 
\begin{align*}
    \frac{1}{\epsilon}&=\mathopen|{c}_{\btheta}'(1)\mathclose| \leq \sum_{i = 1}^{d_1} \mathopen|W_2^i W_1^i\mathclose|\sech^2( W_1^i + b_1^i)\\
    &\leq \sum_{i = 1}^{d_1} \mathopen|W_2^i W_1^i\mathclose| = \mathopen|W_1\mathclose|\mathopen|W_2\mathclose|,\\
    % \frac{1}{\epsilon^2}& = | {c}_{\btheta}''(x=1)|\\
    % &\leq  2\sum_{i = 1}^{d_1} |W_2^iW_1^i||W_1^i||\tanh( W_1^i + b_1^i)|\sech^2( W_1^i + b_1^i)\\
    % &<2\sum_{i = 1}^{d_1} |W_2^iW_1^i||W_1^i| (OMITTABLE),
\end{align*}
by triangle inequality. Thus, $\frac{1}{\epsilon}$ is bounded by the product of the absolute values of the entries in the weights. When the entries in the weight matrices is large, it is easier for the network to follow \eqref{eq:ration}. Note that the result is similar for the transformed case. For a network of one hidden layer, a larger width means it is easier for $\mathopen|W_1\mathclose|\mathopen|W_2\mathclose|$ to be of order $O(\frac{1}{\epsilon})$ as no single entry in $W_1$ and $W_2$ needs to be of order $O(\frac{1}{\epsilon})$, i.e., it is easier for a network of one hidden layer and larger width to satisfy \eqref{eq:1dcd2}. 
% }

When the initialization method of weights in layers follow the Xavier initialization, entries in weights follow a normal distribution with mean zero. Thus, it is easier for networks with two hidden layers with sufficient width to satisfy ${c}_{\btheta}'(1) = -\frac{1}{\epsilon}$. We observe this phenomenon when we compare the results between networks with one hidden layer of width 2 and networks with one hidden layer of width 100. For the 10 random seeds of initialization, two of the one-layered networks were stuck in a local minimum that resulted in a line output in Figure \ref{fig:seeds}, while none of the two-layered networks were stuck in the linear local minimum.

The analysis quickly becomes complicated for more hidden layers, but we observe the same trend. For example, suppose $c_{\btheta}$ is a network with two hidden layers, i.e., for $x\in [0,1]$, 
$${c}_{\btheta} = \sum_{j}^{d_2}W_3^j\tanh\left(\sum_{i}^{d_1}W_2^{i,j}\tanh(xW_1^i+b_1^i)+b_2^{j}\right)+b_3.$$
Taking derivatives gives
\begin{align*}
{c}_{\btheta}' &= \sum_{j}^{d_2}W_3^j\sech^2\left(\sum_{i}^{d_1}W_2^{i,j}\tanh(xW_1^i+b_1^i)+b_2^{j}\right)\left(\sum_{i}^{d_1}W_2^{i,j}W_1^i\sech^2(xW_1^i+b_1^i)\right).\\
% {c}_{\btheta}^{''} &= -2\sum_{j}^{d_2} W_3^j \sech^2(\sum_{i}^{d_1}W_2^{i,j}\tanh(xW_1^i+b_1^i)+b_2^{j})\tanh(\sum_{i}^{d_1}W_2^{i,j}\tanh(xW_1^i+b_1^i)+b_2^{j})\\
% &(\sum_{i}^{d_1}W_2^{i,j}W_1^i\sech^2(xW_1^i+b_1^i))^2+W_3^j\sech^2(\sum_{i}^{d_1}W_2^{i,j}\tanh(xW_1^i+b_1^i)+b_2^{j})\\
% &(\sum_{i}^{d_1}W_2^{i,j}(W_1^i)^2\sech^2(xW_1^i+b_1^i)\tanh(xW_1^i+b_1^i)).
\end{align*}
Thus,
\begin{align*}
   \frac{1}{\epsilon} &= \mathopen|c_{\btheta}'(1)\mathclose|\leq \sum_{j}^{d_2}\mathopen|W_3^j\mathclose|\sech^2\left(\sum_{i}^{d_1}W_2^{i,j}\tanh(W_1^i+b_1^i)+b_2^{j}\right)\left(\sum_{i}^{d_1}\mathopen|W_2^{i,j}W_1^i\mathclose|\sech^2(W_1^i+b_1^i)\right)\\
    &\leq \sum_{j}^{d_2}\mathopen|W_3^j\mathclose|\left(\sum_{i}^{d_1}\mathopen|W_2^{i,j}W_1^i\mathclose|\right)= \mathopen|W_1\mathclose|\mathopen|W_2\mathclose|\mathopen|W_3\mathclose|.
\end{align*}

Instead of depending only on the weights in the first layer and output layer to ensure that the network fulfills the criteria $c_{\btheta}'(1)=-\frac{1}{\epsilon}$, $\frac{1}{\epsilon}$ depends on the sum of products of entries in the weight matrices. 
\textcolor{black}{
The analysis for more hidden layers is similar. One can find that $\frac{1}{\epsilon}\leq \prod_{i = 1}^L |W_i|$, where $L$ is the number of hidden layers. 
}

% \textcolor{black}{ (moved)}
% The loss function is 
% \begin{align*}
% \mathcal L(\btheta) &= \mathcal L_{u}(\btheta)+ \mathcal L_{r}(\btheta)\\
% & = \frac{1}{2}[{c}_{\btheta}(0)^2+({c}_{\btheta}(1)+1)^2]+\frac{1}{N}\sum_{k=1}^N (-\epsilon {c}_{\btheta}^{''}(x_k)+{c}_{\btheta}'(x_k))^2\\
% & = \frac{1}{2}[{c}_{\btheta}(0)^2+({c}_{\btheta}(1)+1)^2]\\
% &+\frac{1}{N}\sum_{k=1}^N (2\epsilon \sum_{j = 1}^{d_1} W_2^j\sum_{i=1}^{d_1}\tanh(x_k W_1^i + b_1^i)\sech^2(x_k W_1^i + b_1^i){(W_1^i)}^2\\
% &+ \sum_{j = 1}^{d_1} W_2^j\sum_{i=1}^{d_1}\sech^2(x_k W_1^i + b_1^i)(W_1^i))^2.
% \end{align*}
% For $c_{\btheta}(x-1)$, the loss function is
% \begin{align*}
% \mathcal L(\btheta) &= \mathcal L_{u}(\btheta)+ \mathcal L_{r}(\btheta)\\
% & = \frac{1}{2}[{c}_{\btheta}(-1)^2+({c}_{\btheta}(0)+1)^2]+\frac{1}{N}\sum_{k=1}^N (-\epsilon {c}_{\btheta}^{''}(x_k-1)+{c}_{\btheta}'(x_k-1))^2\\
% & = \frac{1}{2}[{c}_{\btheta}(-1)^2+({c}_{\btheta}(0)+1)^2]\\
% &+\frac{1}{N}\sum_{k=1}^N (2\epsilon \sum_{j = 1}^{d_1} W_2^j\sum_{i=1}^{d_1}\tanh(x_k W_1^i + b_1^i-W_1^i)\sech^2(x_k W_1^i + b_1^i-W_1^i){(W_1^i)}^2\\
% &+ \sum_{j = 1}^{d_1} W_2^j\sum_{i=1}^{d_1}\sech^2(x_k W_1^i + b_1^i-W_1^i)(W_1^i))^2.
% \end{align*}
% \textcolor{black}{(maybe take out the loss functions and add the equation below)}
Now that we understand the character of the weight matrices in the network, we further explore the effects of shifting inputs. Some of these effects have already been seen in Figure \ref{fig:firstlayer}. Consider the transformed network with one hidden layer, for $x\in [0,1]$, 
\begin{align*}
c_{\btheta}(T(x)) &= \sum_{i = 1}^{d_1} W_2^i\tanh((x-1) W_1^i + b_1^i)+b_2\\
    & = \sum_{i = 1}^{d_1} W_2^i\tanh(x W_1^i + (b_1^i-W_1^i))+b_2
\end{align*}
Note that the transformation $T(x) = x-1$ leads to a new form of bias in the first layer, $b_1^i-W_1^i$. When $W_1^i$ has a large magnitude and $b_1^i$ does not, $b_1^i-W_1^i$ will be large in magnitude and cause $x_k W_1^i + b_1^i-W_1^i$ to be near zero when $x_k=1$. We observe in Figure \ref{fig:afteraffine} that, for both the transformed and untransformed inputs, at least one entry in $W_1$ is of a large magnitude. Naturally, $xW_1$ forms a steep line that approaches zero near $x=0$, whereas $(x-1)W_1$ forms a steep line that approaches zero near $x=1$. When entries in $b_1$ are small, after applying the activation function, a steep gradient forms near where the line approaches zero. Thus, the input transformation $T(x) = x-1$ induces a tendency for the network approximations to contain a steep gradient near $x=1$, aligning with the analytic solution. The analytic form of the loss functions explains the different locations of the steep boundary layers attained by the networks for transformed and untransformed inputs. However, the large magnitude entries in $W_1$ still cause the training process to be difficult. We next explore conditions to alleviate this difficulty.

We explored the effect of network sizes and shifting inputs on the accuracy of approximations produced by the network so far. But we do not need to restrict the transformations on inputs to only shifting. We explore the effect of linear transformation, including stretching, on the performance of PINNs next. 

\begin{figure}[t]
\centering
\begin{subfigure}{0.48\textwidth}
\centering
\includegraphics[width = \linewidth]{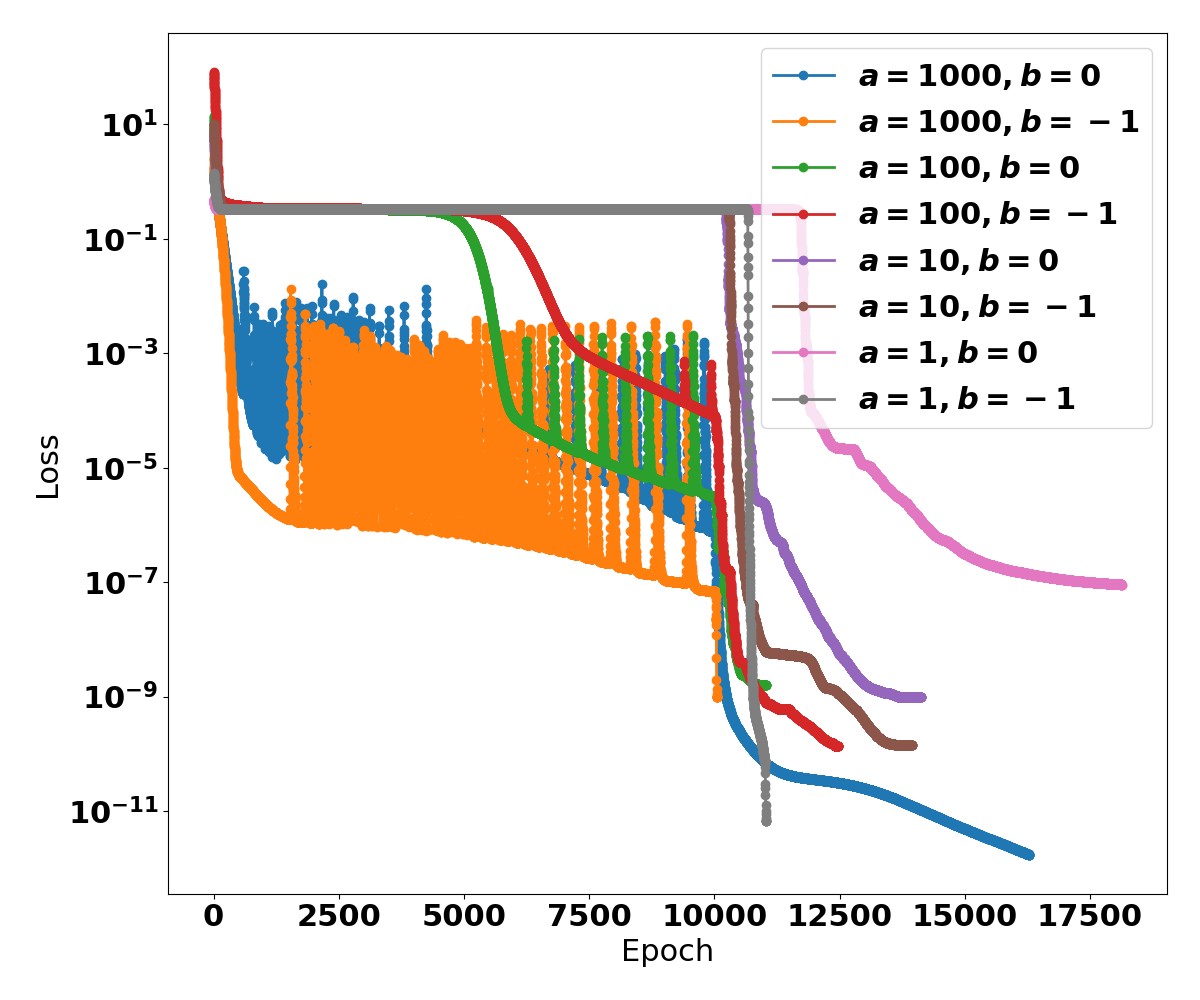}
\caption{}
\label{fig:loss_ab}
\end{subfigure}
\begin{subfigure}{0.48\textwidth}
\centering
\includegraphics[width = \linewidth]{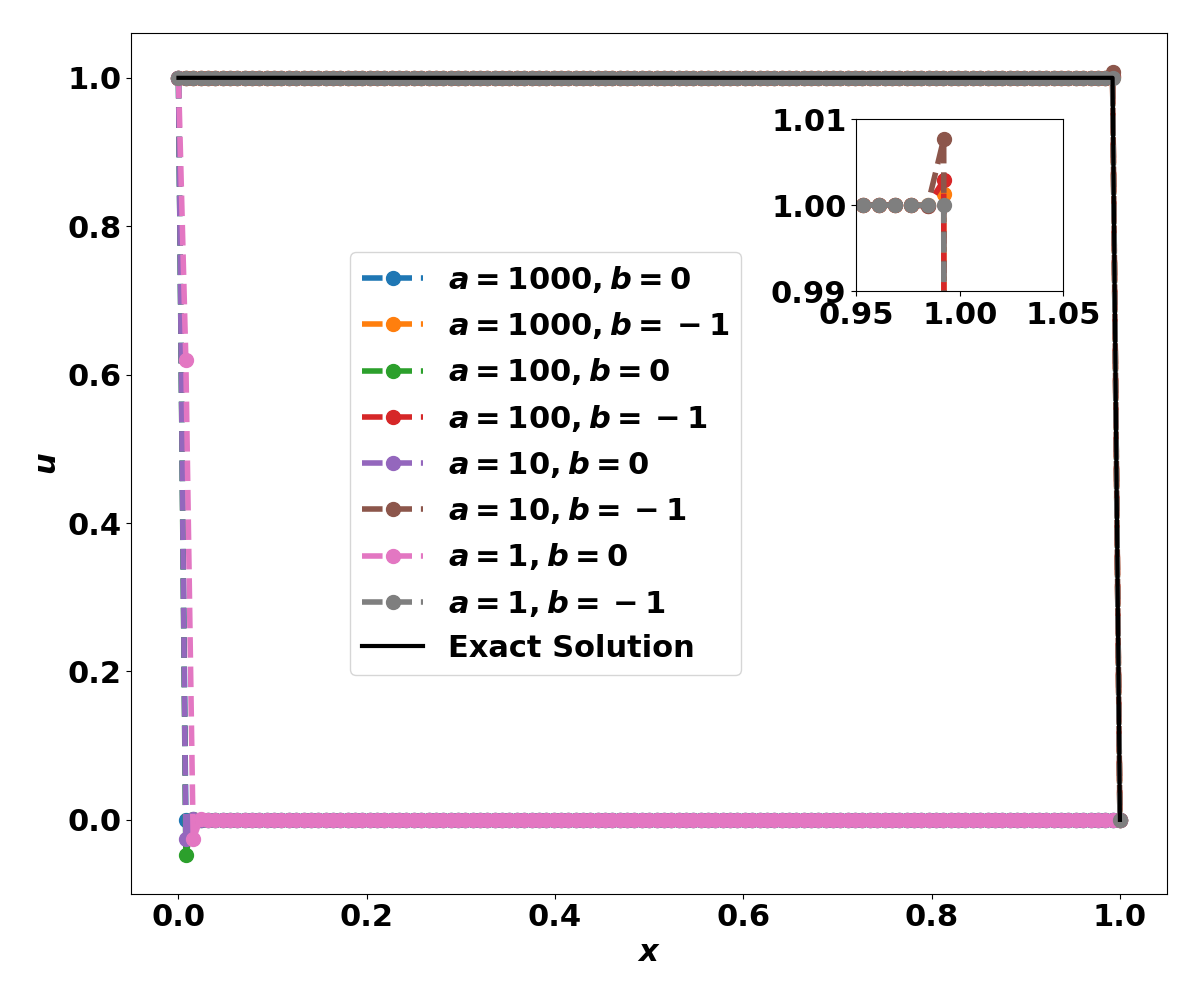}
\caption{}
\label{fig:loss_ap}
\end{subfigure}
\caption{Loss value during training for a network of one hidden layers of width 20 in Figure \ref{fig:loss_ab}. Corresponding approximations in Figure \ref{fig:loss_ap}. First 10000 epochs used Adam optimizer, followed by LBFGS optimizer. Here $a$ and $b$ are the linear transformation parameters in \eqref{eq:u_linear}. }
\end{figure}
\subsection{Input Transformation: Shifting and Scaling}\label{sec:input_transformation}
% % Besides employing a deeper and wider network to fulfill the residual function, we could manipulate the inputs of the network. We propose an input transformation scheme. For spatial coordinates $x$, we perform a transformation of $T(x)$on the input coordinates. Here we will use a linear transformation, i.e. $T(x) = a(x+b)$. 

% % Let the new network structure be
% % \begin{equation}
% % {c}_{\btheta} = W_{L+1}\sigma(W_L \sigma(W_{L-1}\cdots \sigma(W_1(a(x+b))+b_1)+\cdots+b_{L-1})+b_{L})+b_{L+1}.
% % \end{equation}
% % where $a,b\in \R$.
In this section, we examine the effects of input transformations more closely. To simplify the analysis, we consider a network of one hidden layer and width one. As in Section \ref{sec:nn}, suppose the input $x$ goes through a linear transformation $T(x) = a(x+b)$, i.e., for $x\in [0,1]$ and $W_1, W_2, b_1, b_2\in \R$, 
\begin{equation}
c_{\btheta}(T(x)) = \tanh\left(a(x+b)W_1+b_1\right)W_2+b_2.
\label{eq:u_linear}
\end{equation}
Then taking derivatives, we get
\begin{align*}
c_{\btheta}'(T(x)) =  aW_1W_2\sech^2(a(x+b)W_1+b_1),\\
c_{\btheta}^{''}(T(x)) =  -2(aW_1)^2W_2\tanh(a(x+b)W_1+b_1)\sech^2(a(x+b)W_1+b_1).
\end{align*}
Again, assuming $W_1,W_2>0$ and $c_{\btheta}'(x)\neq 0$, we have
$$\frac{c_{\btheta}''(T(x))}{c_{\btheta}^{'}(T(x))}\leq \frac{4}{3\sqrt{3}}aW_1.$$
Instead of only depending on the weights, we could set $a$ to be large, like $\frac{1}{\epsilon}$, to enforce $\frac{c_{\btheta}''(T(x))}{c_{\btheta}^{'}(T(x))} = \frac{1}{\epsilon}$. \\
% explain more, refer to figure`
We also see success experimentally. Figure \ref{fig:loss_ab} shows the loss value during training for different value of $a$. It is clear that as $a$ increases, loss values decrease more rapidly during Adam optimization. Networks also were trained to smaller loss values with larger values of $a$ such as 100 and 1000. Note that even though the training appears to be more successful with larger $a$, regardless of the value of $b$, networks only attain the correct approximations when $b = -1$. 
% Side note, with $a = 1/\epsilon$, $b = -1$, the network of one hidden layer of width 1 can accurately approximate the solution. \\

To understand why different values of transformation bias $b$ lead to different surrogate functions, we will apply the recently developed Neural Tangent Kernel (NTK) theory \cite{jacot2018neural} to PINN \cite{wang2022and}. We follow the work developed in \cite{wang2022and} in the following discussion. 

In \cite{wang2022and}, minimizing the loss function by gradient descent with an infinitesimally small learning rate yields the continuous time gradient flow system
\begin{equation}
\frac{d\theta}{dt} = -\nabla \mathcal{L} (\btheta).
\label{eq:gradientflow}
\end{equation}
%explain the training in gradient vased method relate to kernel regression

% Since training using gradient-based methods is the same as performing kernel regression, information on kernels can serve as an indicator of how neural networks perform as surrogate functions. 
\begin{figure}[t]
\centering
\begin{subfigure}{0.48\textwidth}
\centering
\includegraphics[width = \linewidth]{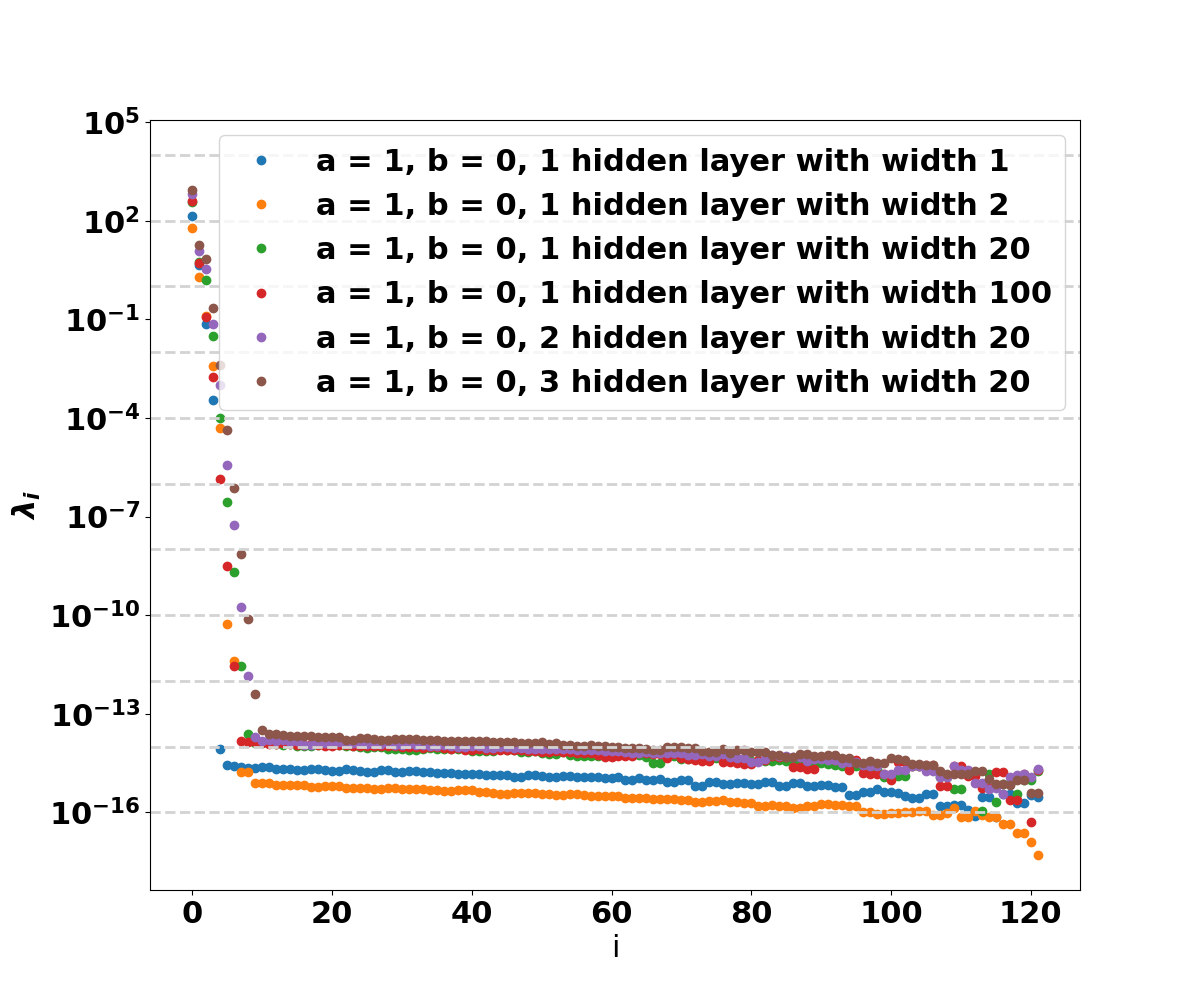}
\caption{}
\label{fig:eigval_width}
\end{subfigure}
\begin{subfigure}{0.48\textwidth}
\centering
\includegraphics[width = \linewidth]{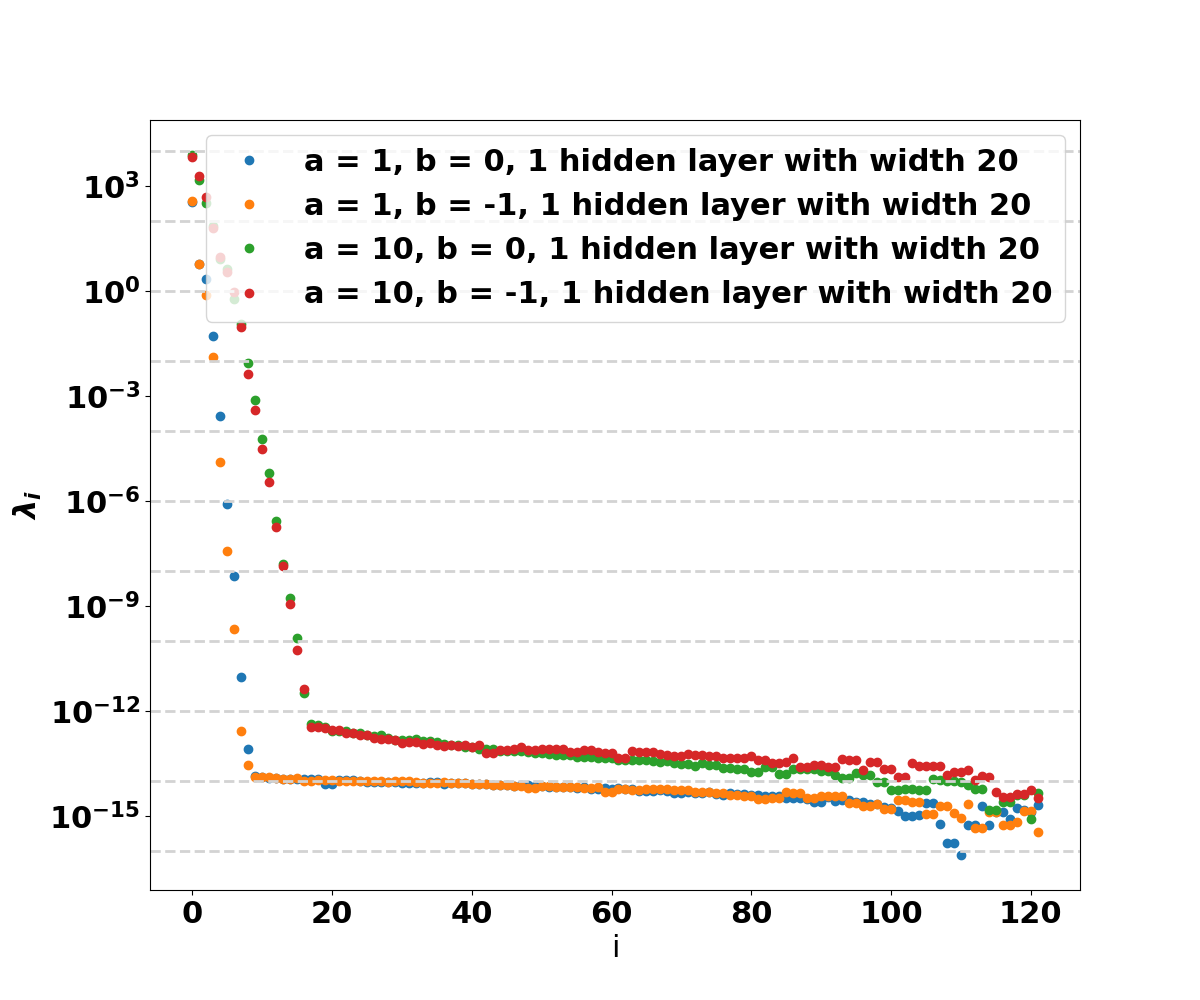}
\caption{}
\label{fig:eigval_trans}
\end{subfigure}
\caption{Eigenvalues of $\bfK$ for different network width and depth (\ref{fig:eigval_width}), and for different transformation parameters $a,b$ (\ref{fig:eigval_trans}). }
\label{fig:eig_val}
\end{figure}

Following the derivation of \cite{jacot2018neural}, the neural tangent kernel can be generally defined at any time $t$ as the neural
network parameters $\theta(t)$ vary during training by gradient descent. The kernel is defined by 
\begin{equation}
{\bfK_{er}}_t(x,x') = \Big\langle\frac{\partial f(x,\btheta(t))}{\partial\btheta},\frac{\partial f(x',\btheta(t))}{\partial\btheta}\Big\rangle.
\end{equation}
% and this kernel converges in probability to a deterministic kernel at random initialization as the width of hidden layers goes to infinity \cite{jacot2018neural}. 
As a consequence of \cite{jacot2018neural}, a properly randomly initialized and sufficiently wide deep neural network trained by gradient descent is equivalent to a kernel regression with a deterministic kernel.

Let $u_s(t) = u_s(x_u, \btheta(t)) = \{u_s(x_u^i,\btheta(t))\}_{i=1}^{N_b}$ and 
% $\mathcal L_r(t) = \mathcal L_r(r_{\btheta(t)}(x_r^i)) =\{\mathcal  L_r(r_{\btheta(t)}(x_r^i))\}_{i=1}^{N_r}$. 
$\mathcal N [u_s(t)] = \mathcal N [u_s(x_r,\btheta(t))] =\{\mathcal N [u_s(x_r^i,\btheta(t))]\}_{i=1}^{N_r}$. 
By \cite{wang2022and}, with gradient flow in \eqref{eq:gradientflow}, $u_s(t)$ and $\mathcal N [u_s(t)]$ obey the following evolution 
\begin{equation}
\begin{bmatrix}
\frac{du_s(x_u, \btheta(t))}{dt}\\
\frac{d\mathcal N u_s(x_r, \btheta(t))}{dt}
\end{bmatrix} = -\begin{bmatrix}
\bfK_{uu}(t) &\bfK_{ur}(t)\\
\bfK_{ru}(t) & \bfK_{rr}(t)
\end{bmatrix}
\begin{bmatrix}
u_s(x_u, \btheta(t))-u_{bc}\\
\mathcal{N}[u_s(x_r,\btheta(t))]-f(x_r)
\end{bmatrix},
\end{equation}
where $\bfK_{ru}(t) = \bfK_{ur}^T(t)$, $\bfK_{uu}(t)\in\R^{N_b\times N_b}$, $\bfK_{ur}(t)\in\R^{N_b\times N_r}$, and $\bfK_{rr}(t)\in\R^{N_r\times N_r}$ are given by 				
\begin{equation}
\begin{split}
(\bfK_{uu})_{ij}(t) &= \Big\langle \frac{du(x_u^i, \btheta(t))}{d\btheta},\frac{du(x_u^j, \btheta(t))}{d\btheta}\Big\rangle \\
(\bfK_{ur})_{ij}(t) &= \Big\langle \frac{du(x_u^i, \btheta(t))}{d\btheta}, \frac{d\mathcal N [u(x_r^j, \btheta(t))]}{d\btheta}\Big\rangle\\
(\bfK_{rr})_{ij}(t) &= \Big\langle \frac{d\mathcal N[u(x_r^i, \btheta(t))]}{d\btheta}, \frac{d\mathcal N[u(x_r^j, \btheta(t))]}{d\btheta}\Big\rangle,\\
\end{split}
\end{equation}
where $\langle \cdot, \cdot \rangle$ denotes the inner product over all trainable parameters in $\btheta$. For example, 
$$(\bfK_{uu})_{ij}(t) = \sum_{\theta\in \btheta} \frac{du(x_u^i, \btheta(t))}{d\theta} \cdot \frac{du(x_u^j, \btheta(t))}{d\theta}. $$

These kernel matrices form a kernel matrix $\bfK$, where 
\begin{equation}
\bfK = \begin{bmatrix}
\bfK_{uu}(t) &\bfK_{ur}(t)\\
\bfK_{ru}(t) & \bfK_{rr}(t)
\end{bmatrix}.
\end{equation}
\begin{figure}[t]
\centering
\begin{subfigure}{0.48\textwidth}
\centering
\includegraphics[width = \linewidth]{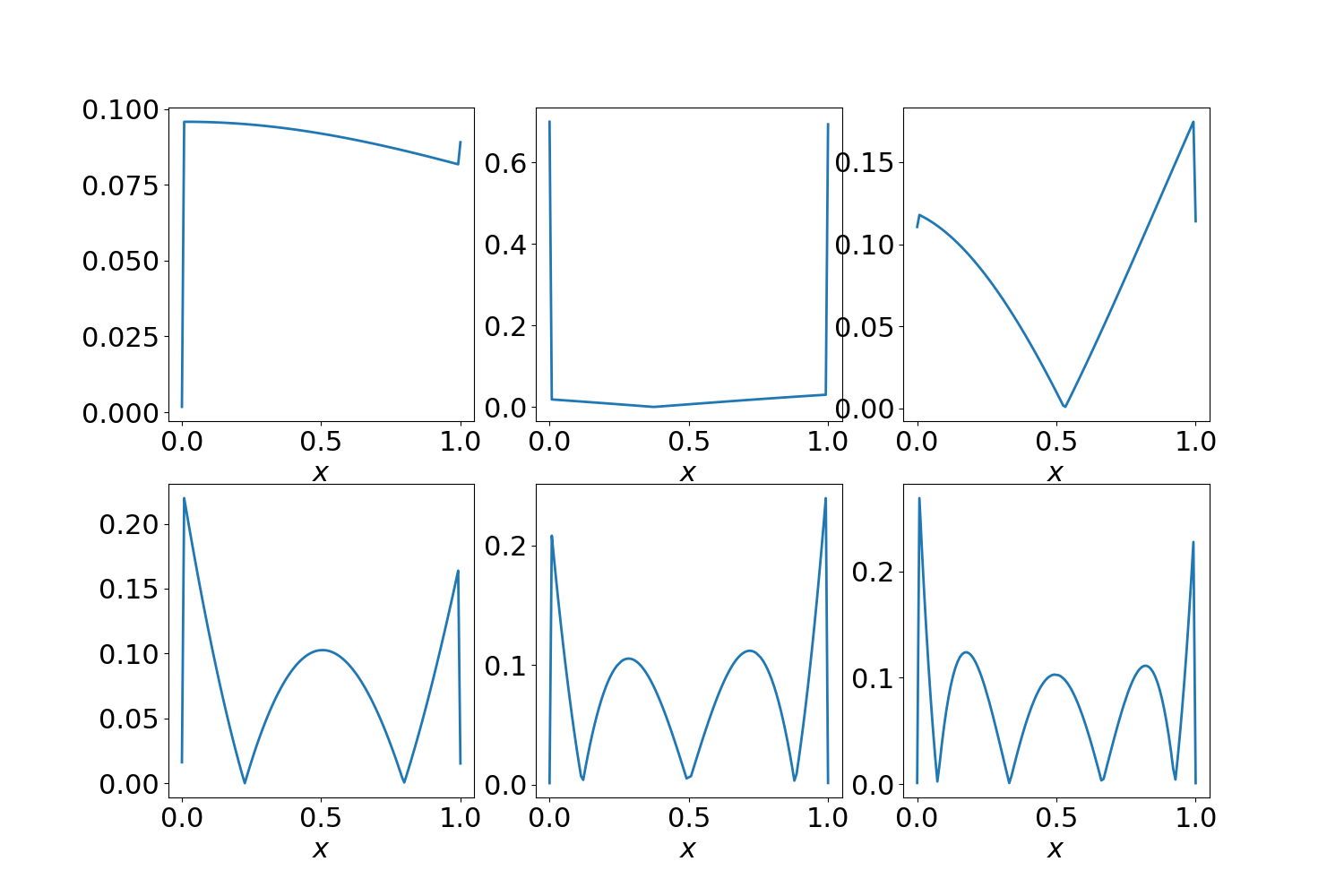}
\caption{$a = 1, b= 0$}
\label{fig:eigvec_none}
\end{subfigure}
\begin{subfigure}{0.48\textwidth}
\centering
\includegraphics[width = \linewidth]{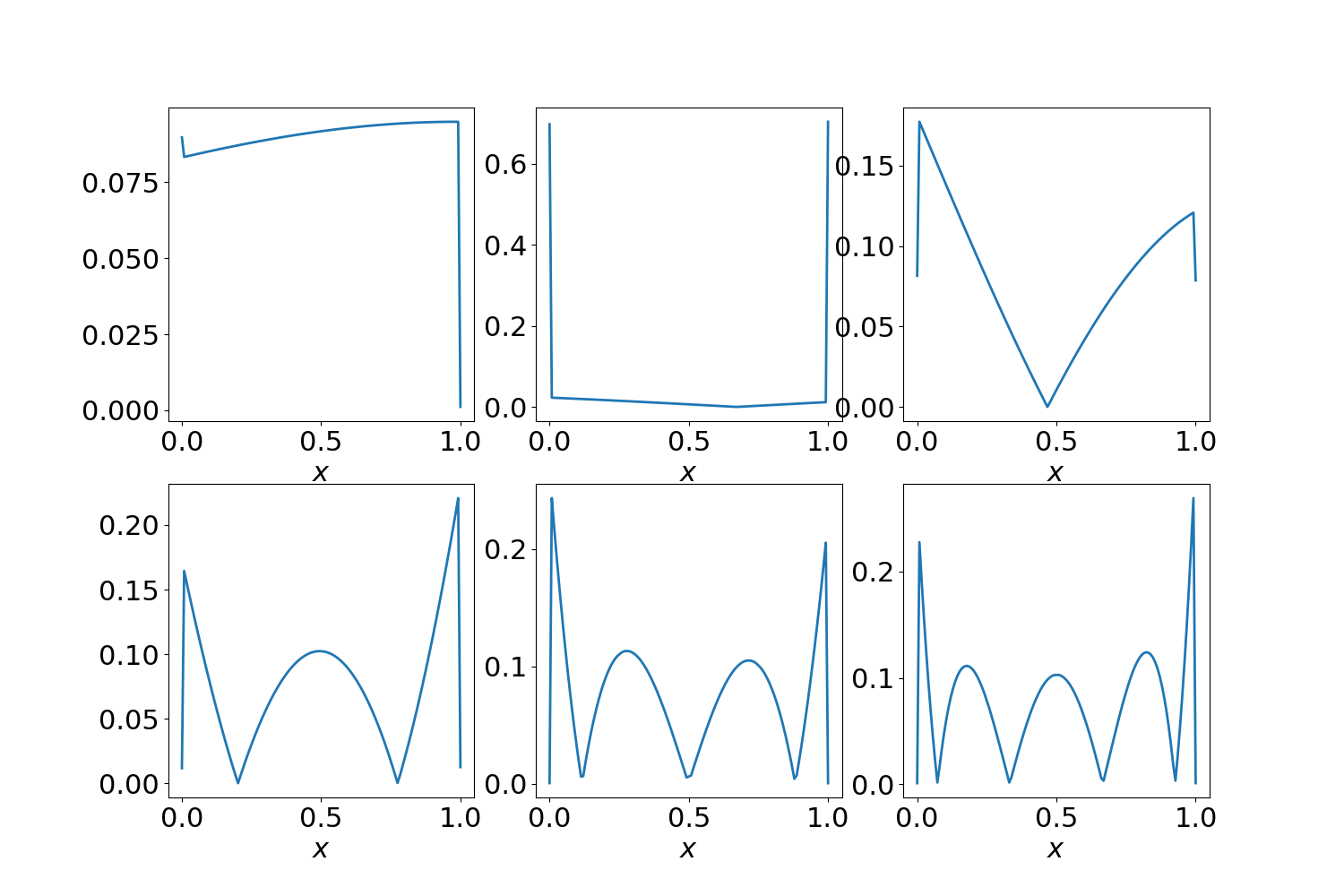}
\caption{$a = 1, b= -1$}
\label{fig:eigvec_inputm1}
\end{subfigure}
\begin{subfigure}{0.48\textwidth}
\centering
\includegraphics[width = \linewidth]{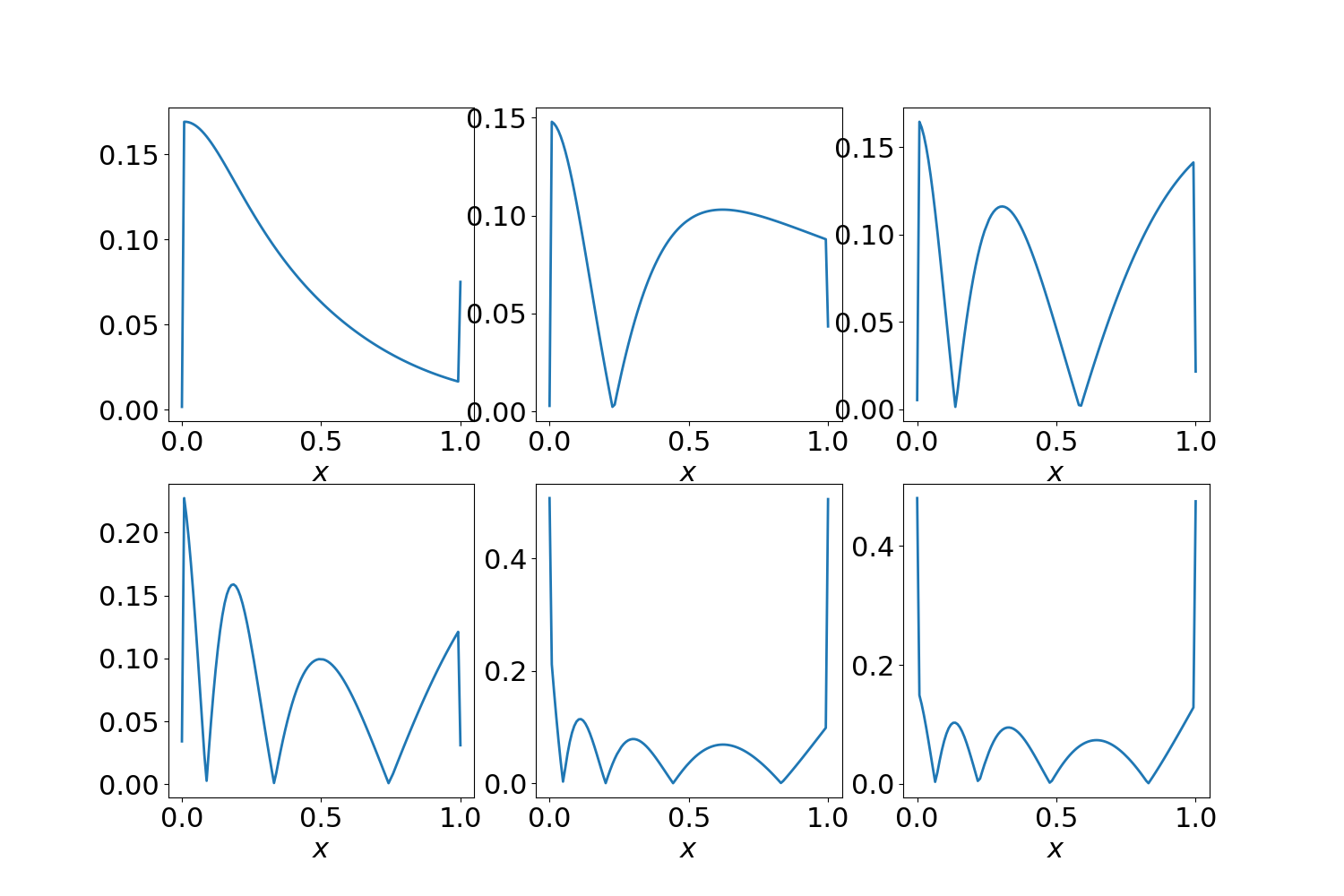}
\caption{$a = 10, b= 0$}
\label{fig:eigvec_input10}
\end{subfigure}
\begin{subfigure}{0.45\textwidth}
\centering
\includegraphics[width = \linewidth]{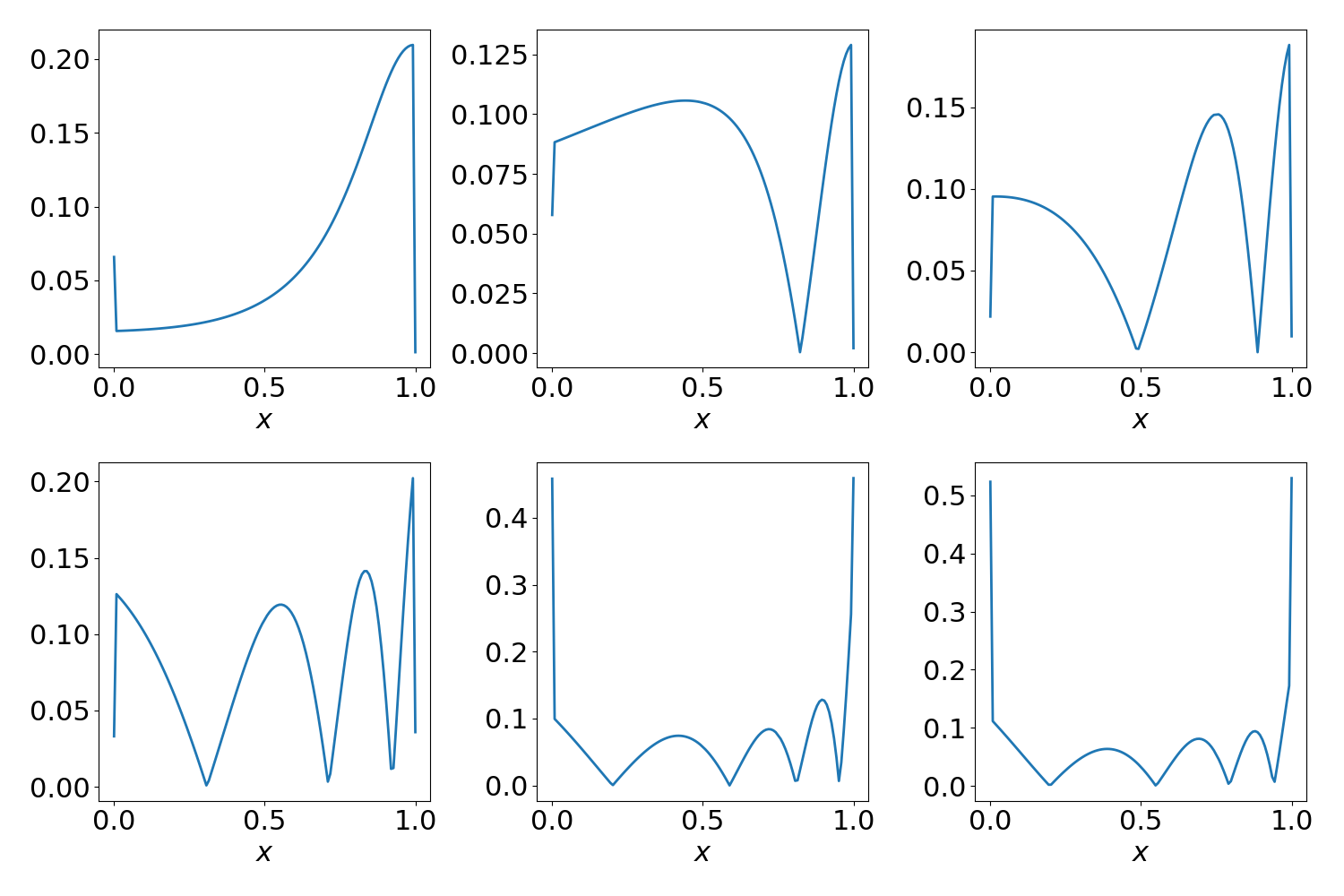}
\caption{$a = 10, b= -1$}
\label{fig:eigvec_input10m10}
\end{subfigure}
\caption{Eigenvectors corresponding to the 6 largest eigenvalues of $\bfK$, for eigenvalues decreasing from left to right, top to bottom. Fig \ref{fig:eigvec_none}: the eigenvectors for $a = 1, b = 0$. Fig \ref{fig:eigvec_inputm1}: the eigenvectors for $a = 1, b = -1$. Fig \ref{fig:eigvec_input10}: the eigenvectors for $a = 10, b = 0$. Fig \ref{fig:eigvec_input10m10}: the eigenvectors for $a = 10, b = -1$}
\label{fig:eig_vec}
\end{figure}
Let $N = N_b+N_r$. By construction, $\bfK$ is a symmetric positive semi-definite matrix of size $N\times N$. In addition, the NTK of PINNs converge in probability to a deterministic limiting kernel when the network width $d\rightarrow \infty$ \cite{wang2022and}. Thus, it is enough to study the behavior of $\bfK$ at initialization. Consider the eigen-decomposition of $\bfK$ at initialization, $\bfK = \bfQ^T\Lambda \bfQ$, where $\bfQ$ is an orthogonal matrix whose columns are eigenvectors of $\bfK$ and $\Lambda$ is a diagonal matrix whose entries are the corresponding eigenvalues $\lambda_i\geq 0$ of $\bfK$. As shown in \cite{wang2022and}, 

\begin{equation}
\bfQ\left(
\begin{bmatrix}
\frac{du(x_u, \btheta(t))}{dt}\\
\frac{d\mathcal N[u(x_r, \btheta(t))]}{dt}
\end{bmatrix}-\begin{bmatrix}
u_{bc}\\
f(x_r)
\end{bmatrix}\right )
\approx -e^{-\Lambda t}\bfQ \cdot
\begin{bmatrix}
u_{bc}\\
f(x_r)
\end{bmatrix}
% \dot 
% \begin{bmatrix}
% g(x_b)\\
% f(x_r)
% \end{bmatrix}
\label{eq:eig_kernel}
\end{equation}

Thus, the $i$-th component of the left-hand side of \eqref{eq:eig_kernel} will decay approximately at the rate of $e^{-\lambda_i t}$. 
\textcolor{black}{Notice that the left hand side is the evolution of the absolute training error mapped by $Q$.  The components of the target function that correspond to kernel eigenvectors with larger eigenvalues will be learned faster. } One can gain insight into the decay process of the training errors by analyzing $\bfK$. % The eigenvector of $\bfK_{rr}$ provides insight into the type of functions to where ${c}_{\btheta}$ will converge. 

% convergence rate definination
From \eqref{eq:eig_kernel}, the training errors decay exponentially at a rate determined by the eigenvalues in $\Lambda$. Let the average convergence rate $c$ be defined as the mean of all the eigenvalues $\lambda_i$, i.e. 
\begin{equation}
c = \frac{\sum_{i=1}^N\lambda_i}{N}= \frac{Tr(\bfK)}{N}.
\end{equation}
\textcolor{black}{The convergence rate measures the how fast the training errors converge to 0.} In particular, \cite{wang2022and} shows that $Tr(\bfK) = Tr(\bfK_{uu})+Tr(\bfK_{rr})$, and \textcolor{black}{it follows that $c = \frac{Tr(\bfK_{uu})+Tr(\bfK_{rr})}{N}$. If $ Tr(\bfK_{rr})\gg Tr(\bfK_{uu})$,  the surrogate neural network converges to the PDE governing equation faster than to the boundary conditions. In practice, without input transformation, we found experimentally that $ Tr(\bfK_{rr})\gg Tr(\bfK_{uu})$ at initialization as $Tr(\bfK_{rr})/Tr(\bfK_{uu})\approx 50$, regardless of the network depth or width.}
%add Krr Kuu eigenvalues
% When applied to the one-dimensional convection-diffusion equation with a network of one hidden layer and width two, we found experimentally that $\bfK_{uu}$ dominated $\bfK_{rr}$ at initialization, where $c_{rr}$ is approximately $0.27$ when $c_{uu}$ is approximately 1.4. However, when the network gets larger, $\bfK_{rr}$ starts to dominate $\bfK_{uu}$. 
% It is great to have NTK as lenses we can look through to examine the effect of input transformations on the behavior of PINNs. 

\begin{figure}[t]
\centering
\begin{subfigure}{0.44\textwidth}
\centering
\includegraphics[width = \linewidth]{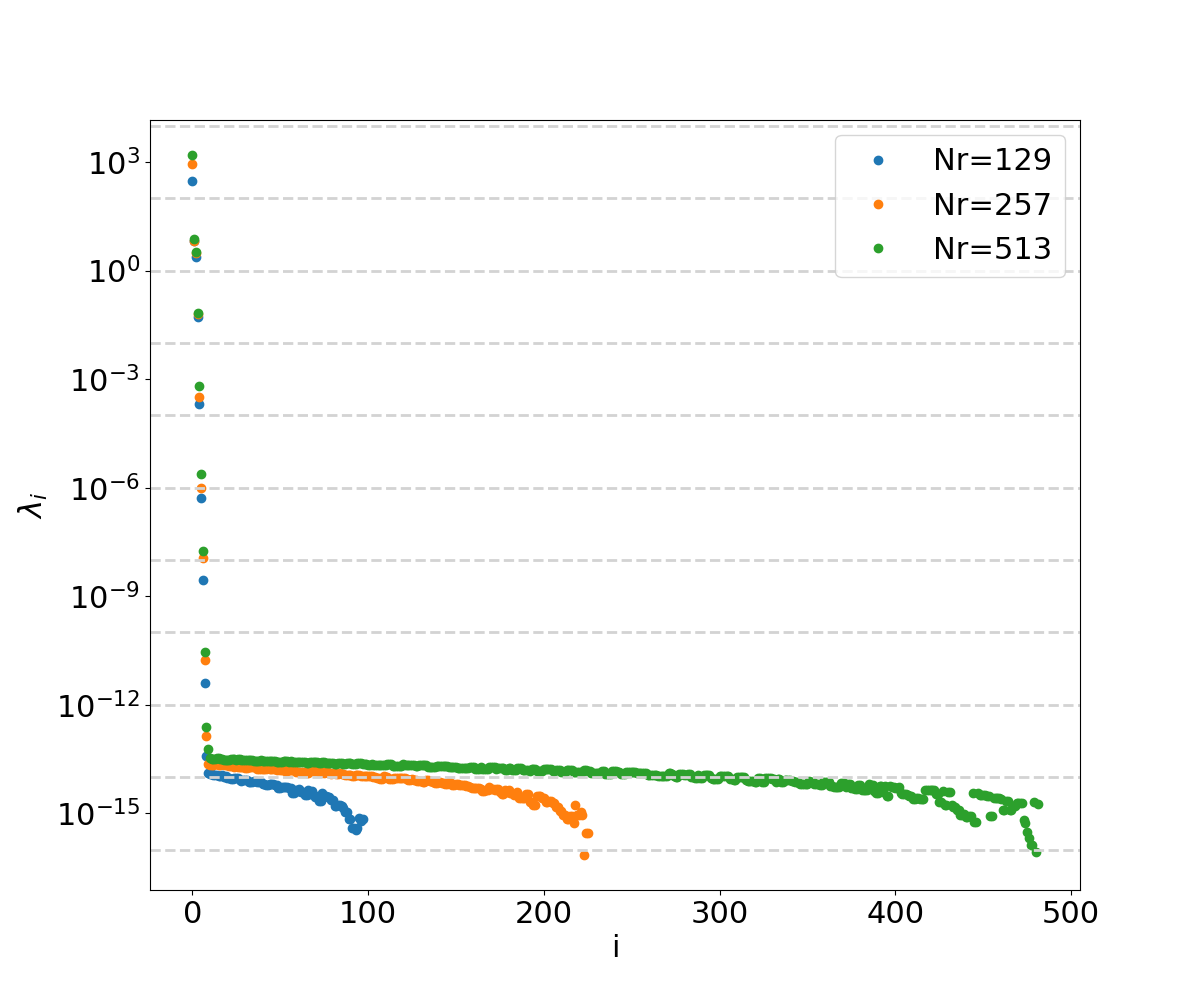}
\caption{}
\label{fig:eigval_Nr}
\end{subfigure}
\begin{subfigure}{0.48\textwidth}
\centering
\includegraphics[width = \linewidth]{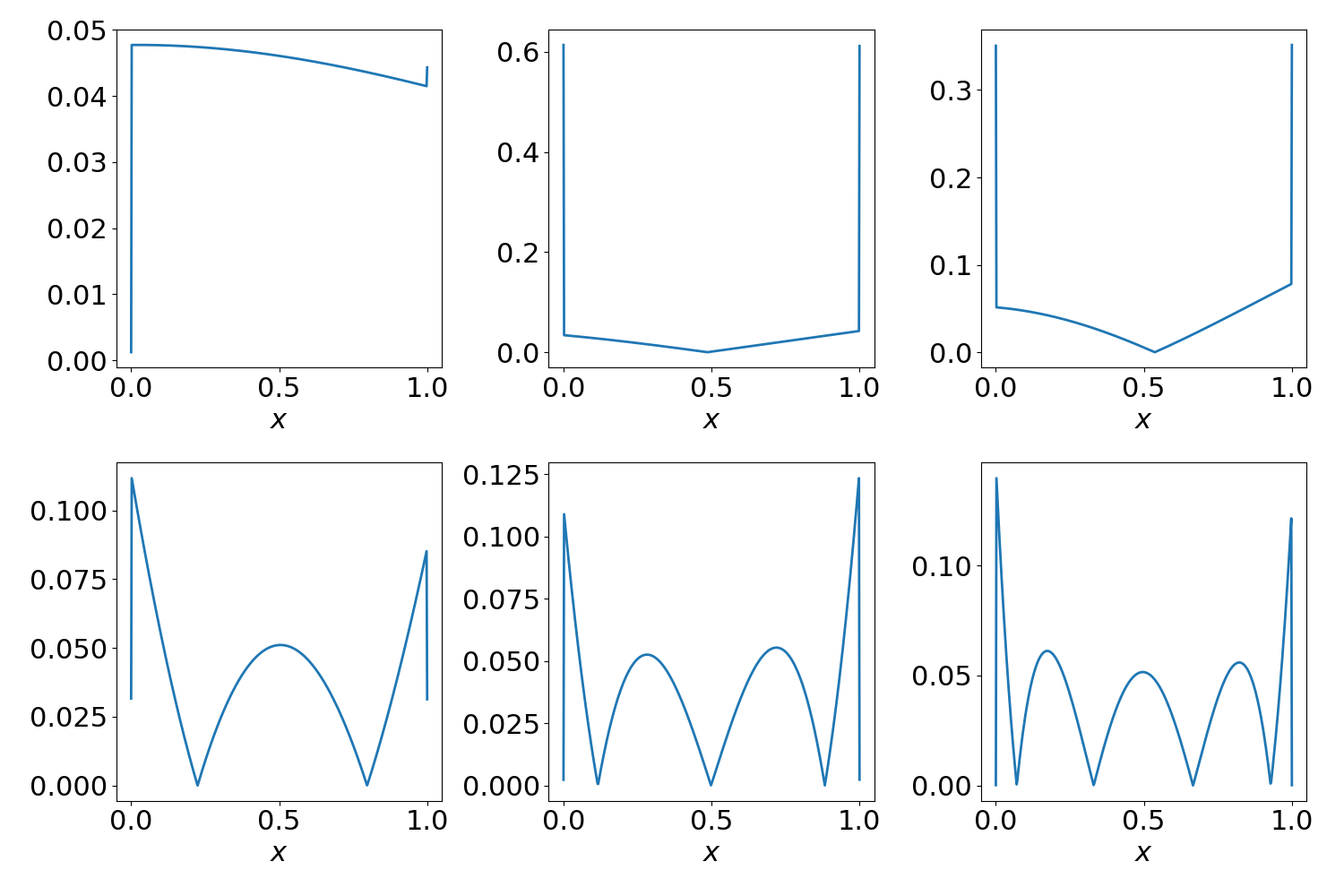}
\caption{}
\label{fig:eigvec_513}
\end{subfigure}
\caption{\ref{fig:eigval_Nr} Eigenvalues of $\bfK$ for different residual sample sizes . \ref{fig:eigvec_513} Eigenvectors of $\bfK$ corresponding to the six largest eigenvalues for sample size 513. }
\label{fig:eig_val_Nr}
\end{figure}

When observing the eigenvalues of $\bfK$, we notice that there are only a few eigenvalues of large magnitude. The situation does not change much when network width or depth changes. Figure \ref{fig:eigval_width} shows that even for deeper and wider networks, the majority of the eigenvalues are close to zero.

% b does not influence the eigenvalue
% a influence
When the linear transformation is applied to the inputs, as shown in Figure \ref{fig:eigval_trans} the eigenvalues of $\bfK$ \textcolor{black}{are similar in magnitude} for different values of $b$. As $a$ increases, the magnitude of eigenvalues increases, which indicates a higher convergence rate. This observation is consistent with the faster convergence of loss values for larger $a$ shown in Figure \ref{fig:loss_ab}. \textcolor{black}{Note that for larger values of $a$, 1000 for example, the loss values oscillate. This is due to a fixed learning rate applied in Adam optimizations. The oscillatory behavior could be reduced with a decaying learning rate. }

The eigenvalues only indicate the convergence rate. But it is still unclear why different values of $b$ lead networks to different approximations. To answer this question, we need to examine the eigenvectors of $\bfK$. Since a majority of the eigenvalues are close to zero, we consider 6 eigenvectors corresponding to the largest 6 eigenvalues. These are plotted in Figure \ref{fig:eig_vec}, which shows a clear distinction between the case of $a = 1, b=0$ and the case of $a = 1, b=-1$. This distinction becomes more dramatic when $a = 10$. Figure \ref{fig:eigvec_input10} suggests that for $b=0$, the eigenvectors corresponding to the six largest eigenvalues themselves tend to have boundary layers near $x=0$, whereas for $b=-1$, these eigenvectors tend to form boundary layers near $x=1$. Since these eigenvectors represent the components in the residual that will decay the fastest, the boundary layers in the top eigenvectors will appear in the approximations obtained by the network. As a result, approximations of PINN will attain boundary layers near $x=0$ for $b=0$ and near $x=1$ for $b=-1$.

To validate the effect of input transformations, we need to be certain that other factors, such as the network size and residual sample size, do not play a vital role. We examine the eigenvalues of $\bfK$ with different residual sample size $N_r = 129,257,513$. As seen in Figure \ref{fig:eig_val_Nr}, the eigenvalues of $\bfK$ do not vary much for different sample sizes and the eigenvectors exhibit similar structures, which translates to the observation that providing more samples does not benefit the training or the accuracy of approximations. 

\begin{figure}[t]
\centering
\begin{subfigure}{0.48\textwidth}
\centering
\includegraphics[width = \linewidth]{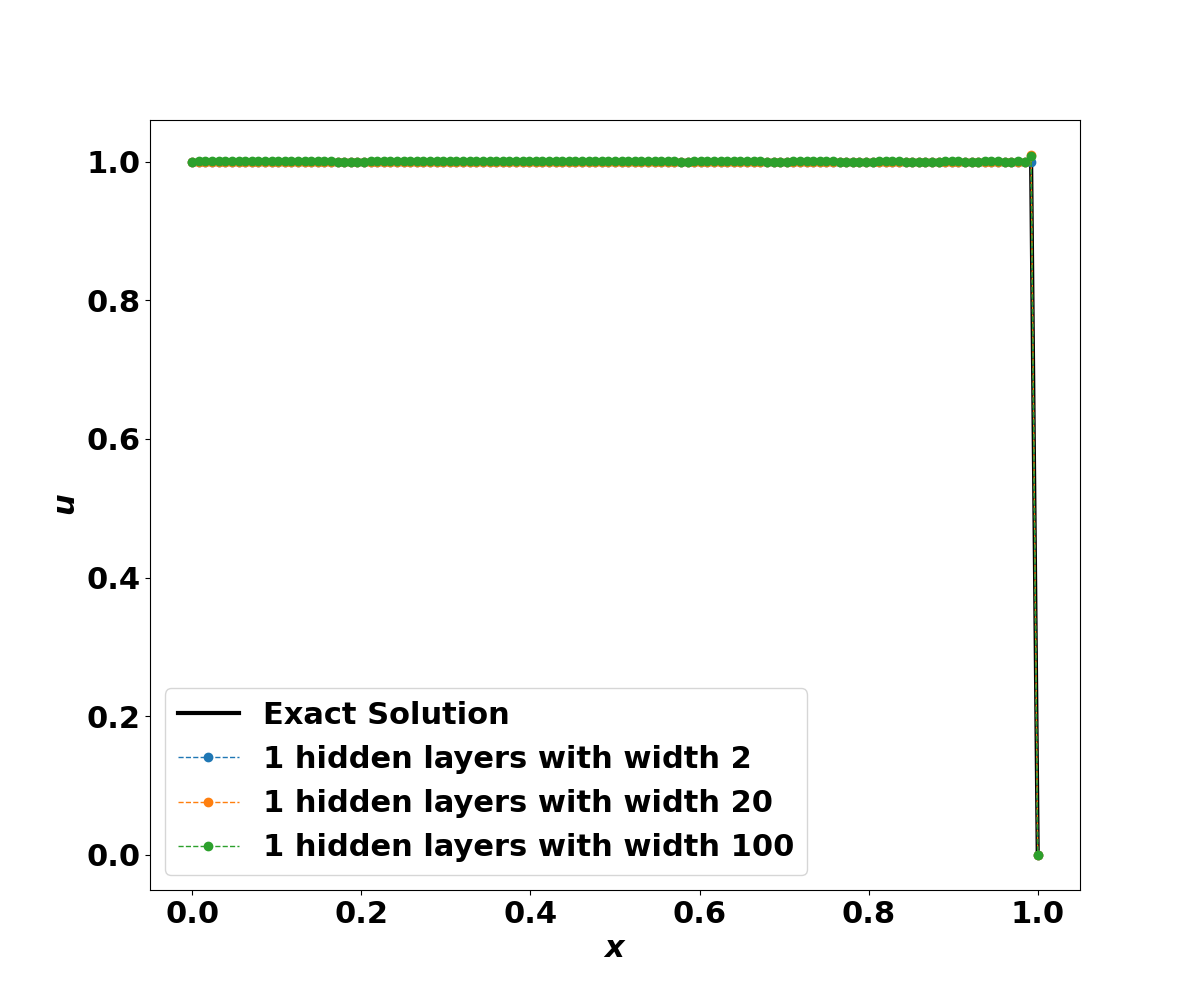}
\caption{}
\label{fig:inputm1_width}
\end{subfigure}
\begin{subfigure}{0.48\textwidth}
\centering
\includegraphics[width = \linewidth]{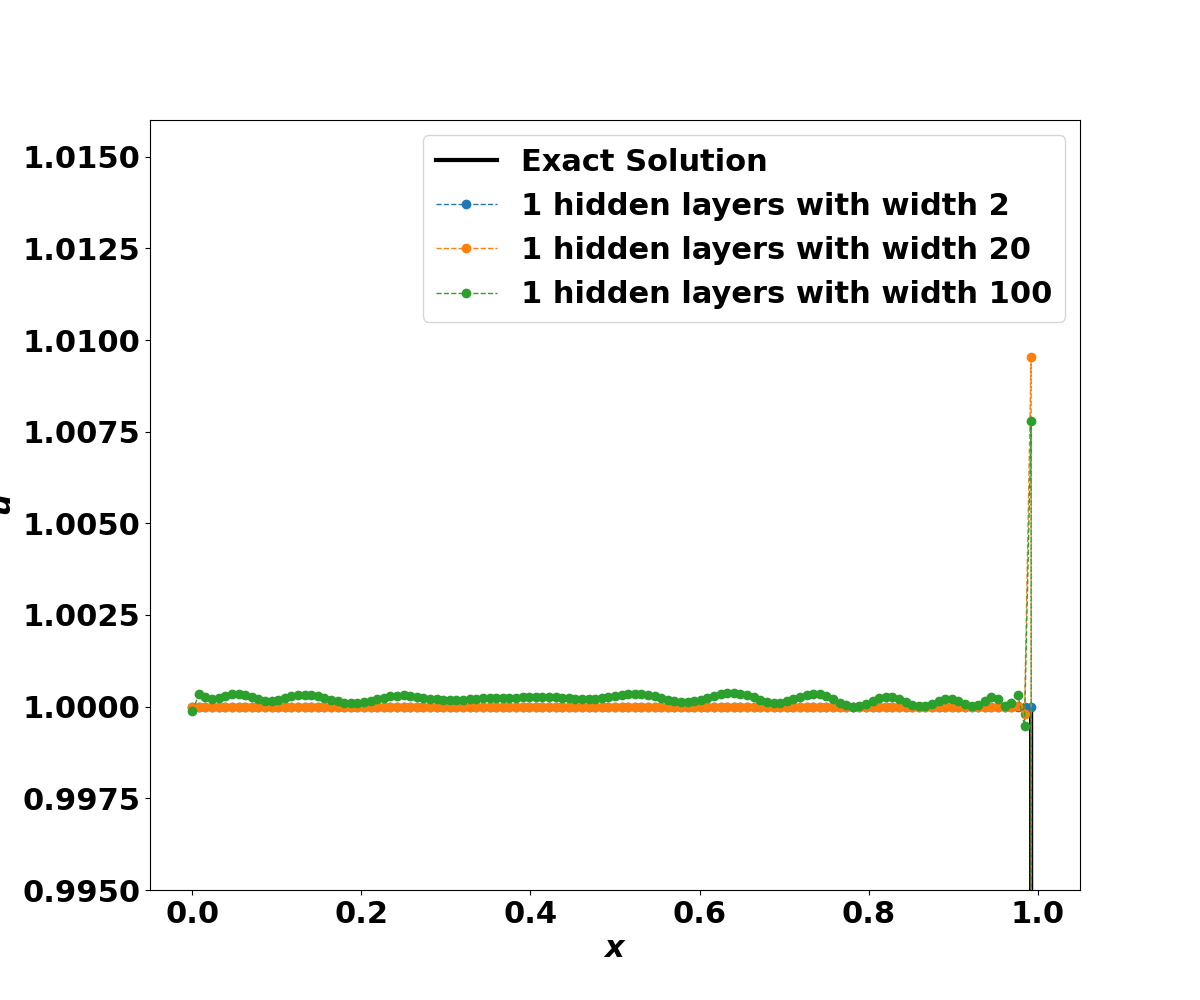}
\caption{}
\label{fig:inputm1_width_zoom}
\end{subfigure}
\caption{\ref{fig:inputm1_width} Approximations produced by trained networks of different size. \ref{fig:inputm1_width_zoom} Zoom in view near $u=1$. }
\label{fig:inputm1_w}
\end{figure}

\begin{figure}[t]
\centering
\begin{subfigure}{0.23\textwidth}
\centering
\includegraphics[width = \linewidth]{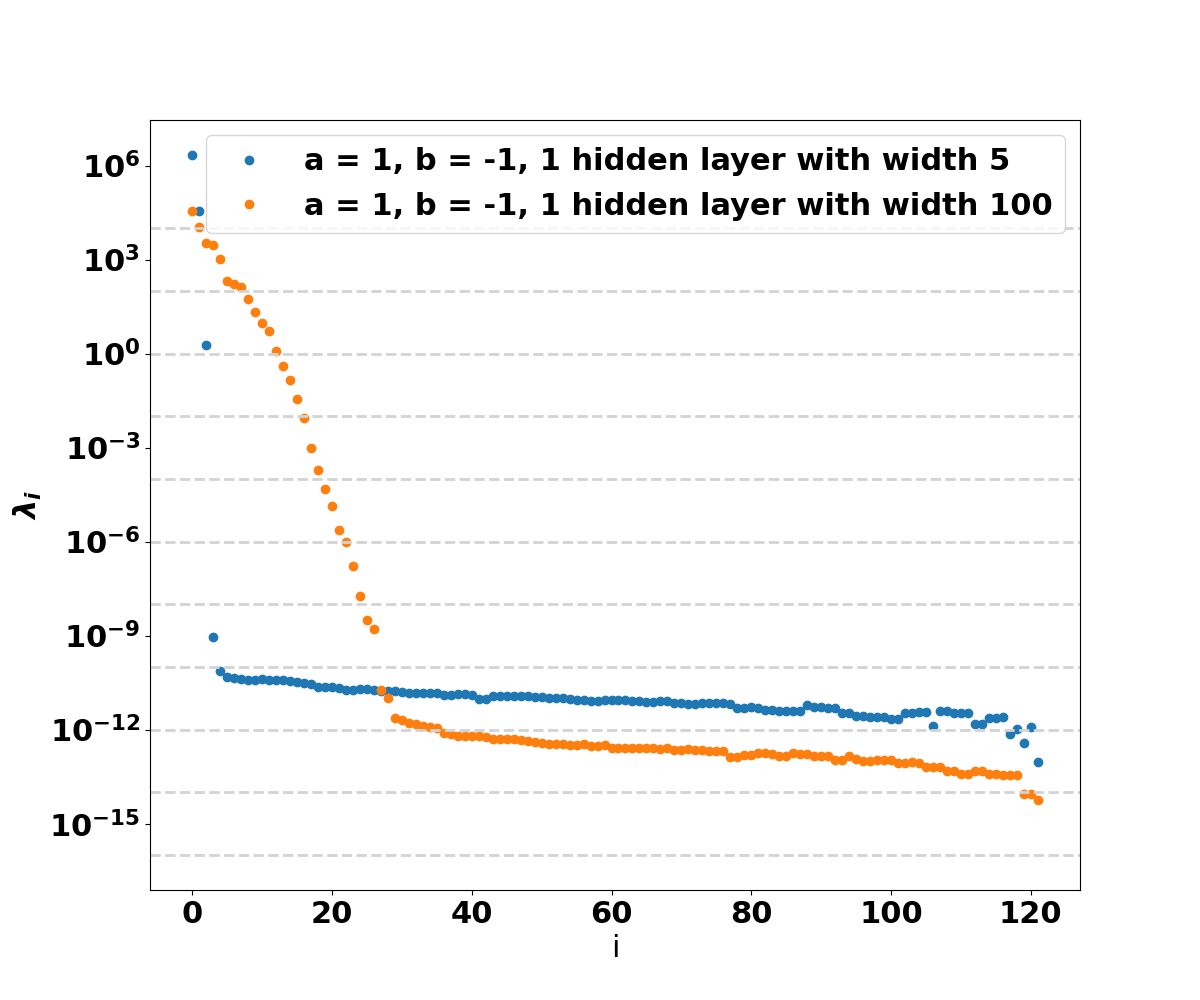}
\caption{}
\label{fig:eigval_width_comp}
\end{subfigure}
\begin{subfigure}{0.38\textwidth}
\centering
\includegraphics[width = \linewidth]{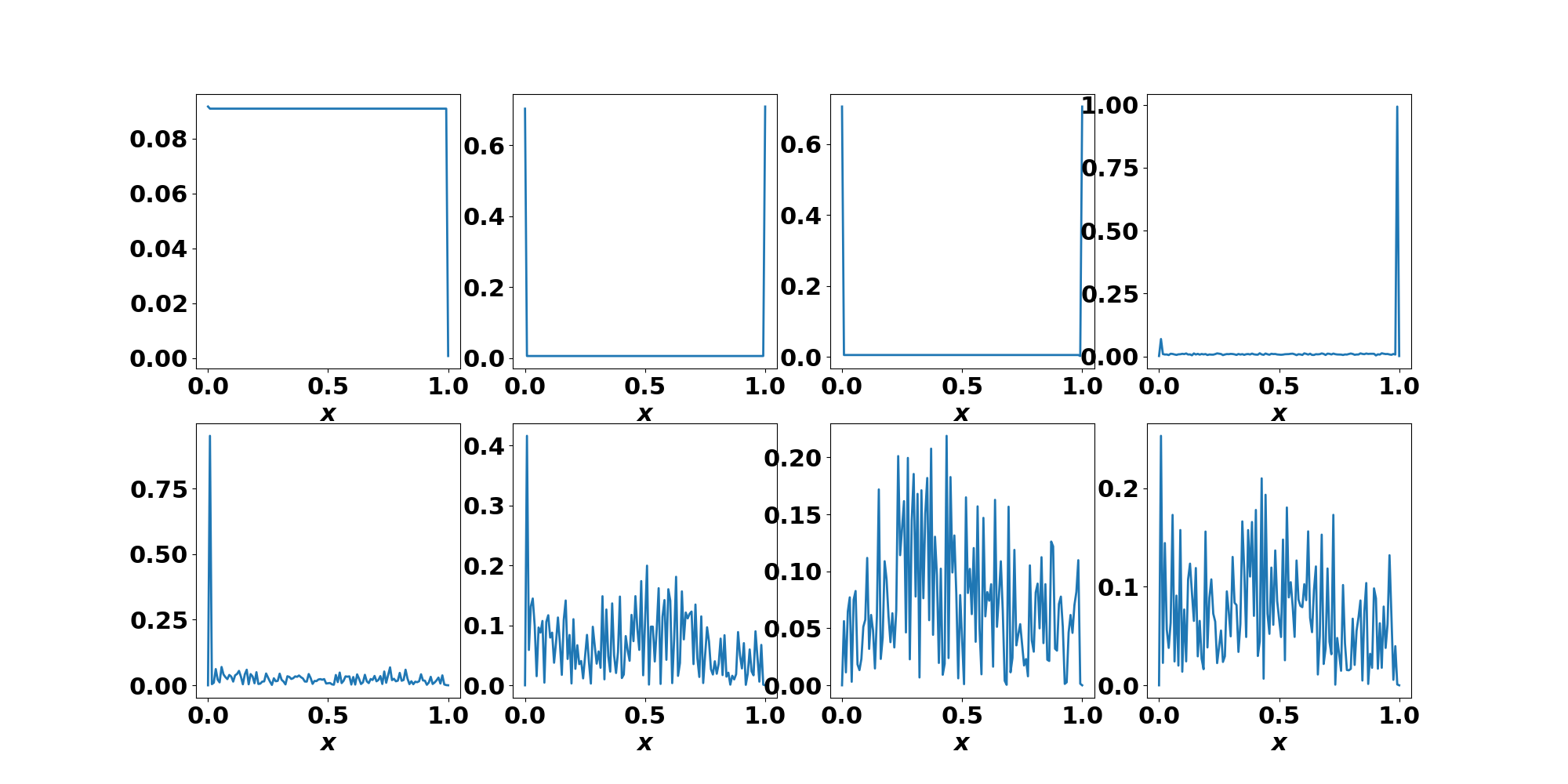}
\caption{}
\label{fig:eigvec_1layer2}
\end{subfigure}
\begin{subfigure}{0.345\textwidth}
\centering
\includegraphics[width = \linewidth]{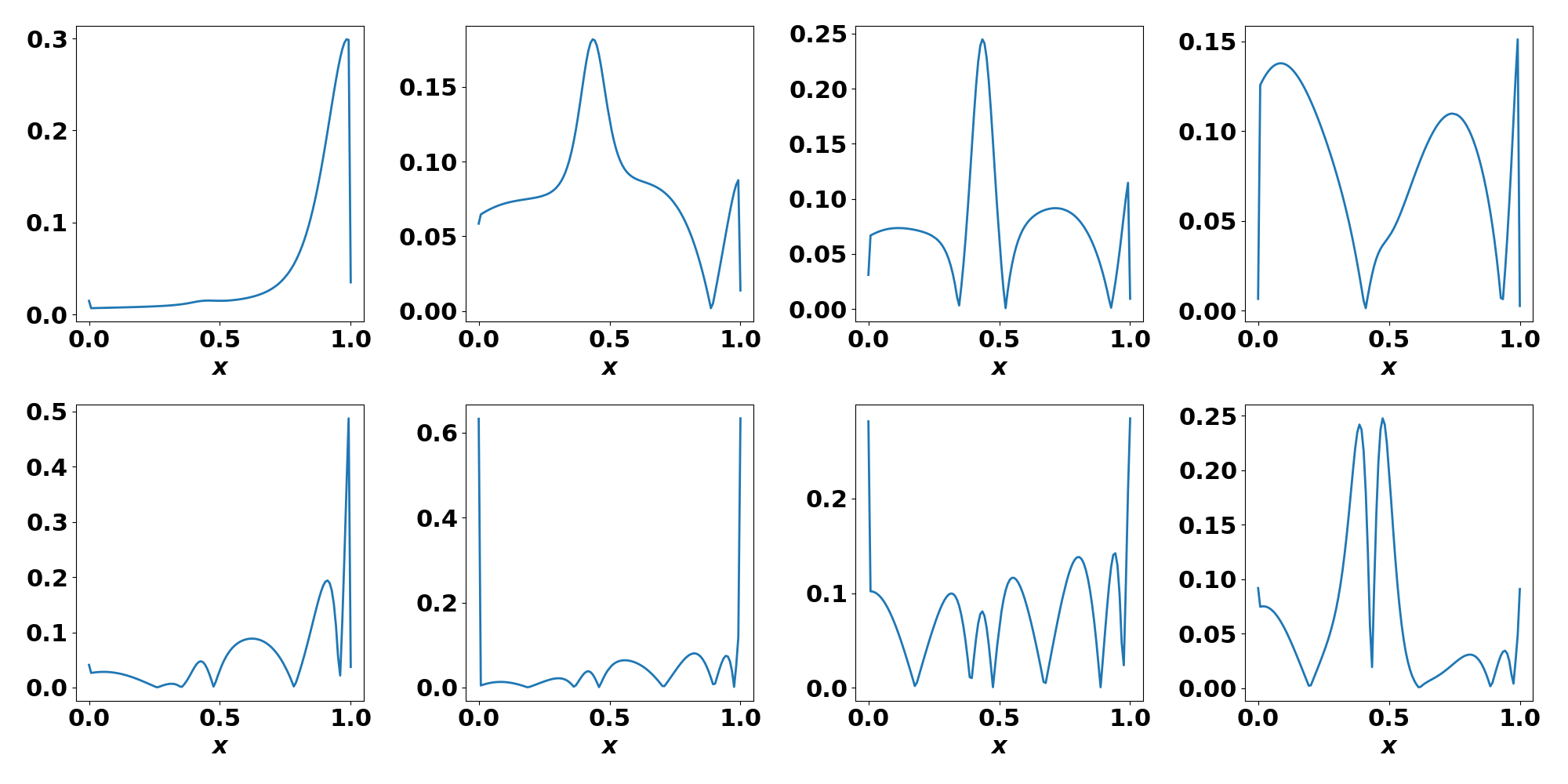}
\caption{}
\label{fig:eigvec_1layer100}
\end{subfigure}
\caption{Eigenvalues and eigenvectors corresponding to the 8 largest eigenvalues of $\bfK$ for $a = 1, b = -1$ after training. Eigenvalues decrease from left to right, top to bottom. Fig \ref{fig:eigval_width_comp}: the eigenvalues of NTK for networks of 1 hidden layer with width 5 and 100, respectively; Fig \ref{fig:eigvec_1layer2}: the eigenvectors for a network of 1 hidden layer with width 5, Fig \ref{fig:eigvec_1layer100}: the eigenvectors for a network of 1 hidden layer with width 100.}
\label{fig:eig_vec_width}
\end{figure}

The width of networks is a tricky subject. Using transformation $x\rightarrow x-1$, we fixed the residual sample size to 129 and trained networks of three different widths. Note that the analytic form yielded that networks with larger width could attain the steep gradient more effectively. But the story took a turn. As shown in Figure \ref{fig:inputm1_w}, the approximations made with network structure of width 2 actually produced the most accurate surrogate exact solution. The network of width 20 started to produce jumps near the steep gradient layer, and the network of width 100 starts to exhibit oscillations throughout the domain. This behavior is caused by the fact that networks of larger width contain higher frequency modes in their eigenvectors of $\bfK$, as shown in Figure \ref{fig:eig_vec_width}. \textcolor{black}{We examine the eigenvectors of $\bfK$ after training to illustrate the different modes the networks converged to. It is clear that the eigenvectors of the network with width 100 display oscillations in the eigenvectors corresponding to the eight largest eigenvalues. The oscillations in the eigenvectors explain the oscillations in the approximations shown in Figure~\ref{fig:inputm1_w}. From Figure \ref{fig:eigval_width_comp}, one can see that only the largest 4 eigenvalues of $\bfK$ for the network of one hidden layer with width 5 are larger than $10^{-9}$, whereas all of the largest 7 eigenvalues of $\bfK$ for the network of width 100 are larger than $10^2$. Thus, the eigenvectors with noisy oscillations corresponding to the fourth to eighth largest eigenvalues in Figure \ref{fig:eigvec_1layer2} are not significant. The oscillations in the eigenvectors in Figure \ref{fig:eigvec_1layer100} are meaningful. }


\section{The Two-Dimensional Convection-Diffusion Equation}
\label{sec:2dcd}
\begin{figure}[t]
\centering
\begin{subfigure}{0.32\textwidth}
\centering
\includegraphics[width =\linewidth]{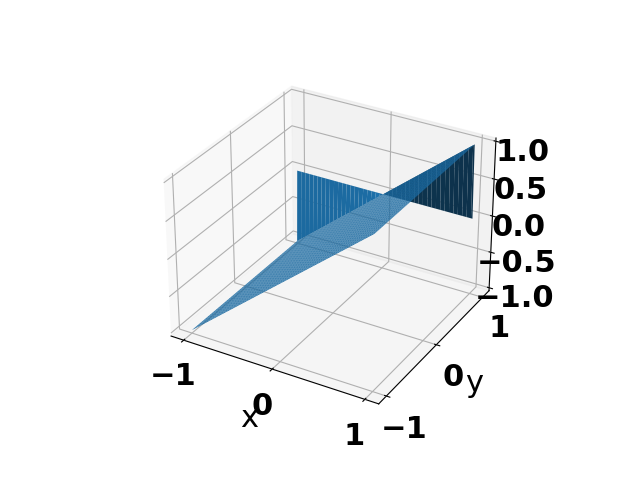}
\caption{}
\label{fig:2dsol}
\end{subfigure}
\begin{subfigure}{0.32\textwidth}
\centering
\includegraphics[width = \linewidth]{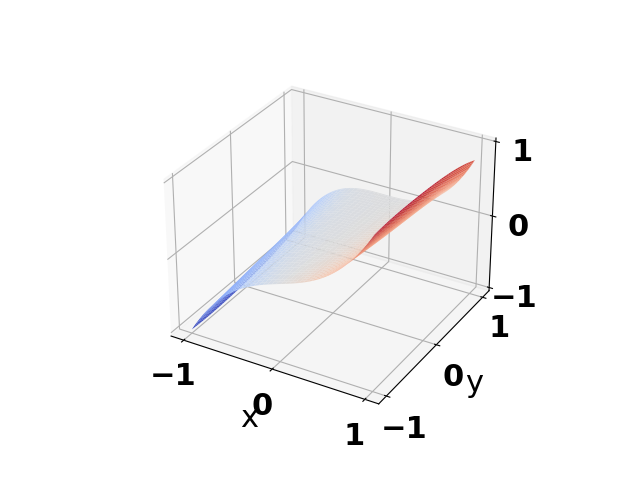}
\caption{}
\label{fig:2d_none}
\end{subfigure}
\begin{subfigure}{0.32\textwidth}
\centering
\includegraphics[width = \linewidth]{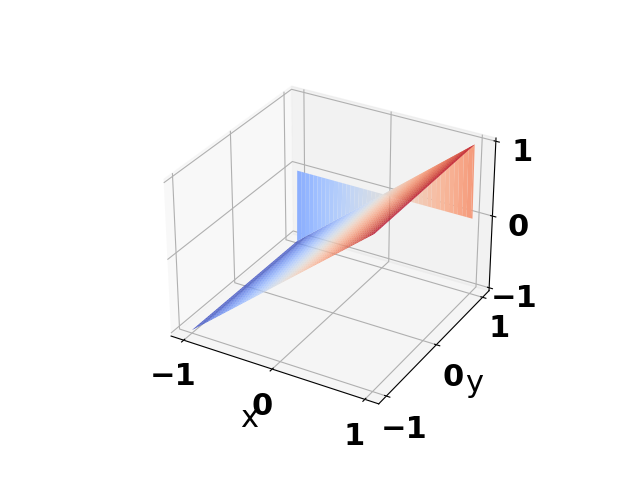}
\caption{}
\label{fig:2d_yinput10m10}
\end{subfigure}
\caption{Solution of \eqref{eq:2d_cd} for $\epsilon = 1e-4$ in \ref{fig:2dsol}. Approximations by PINN in \ref{fig:2d_none} and \ref{fig:2d_yinput10m10}. The input transformations are \ref{fig:2d_none}: no transformation; \ref{fig:2d_yinput10m10}: $(x,y)\rightarrow (x,10(y-1))$.}
\label{fig:2d_pinn}
\end{figure}

We can extend our understanding of the performance of PINN to a simple two-dimensional case, where there is only one steep gradient boundary layer in the domain. Consider the following equation:
\begin{equation}
\begin{split}
-\epsilon \Delta u+\vec{\omega}\cdot \nabla u = 0, \quad \text{ for }(x,y)\in (-1,1)\times (-1,1)\\
u(x, -1) = x,\quad u(x,1) = 0,\\
u(-1,y) \approx -1,\quad u(1,y) \approx 1,
\end{split}
\label{eq:2d_cd}
\end{equation}
where $\vec{\omega} = (0,1)$. The exact solution of \eqref{eq:2d_cd} is 
$$u(x,y) = x\left(\frac{1-e^{(y-1)/\epsilon}}{1-e^{-2/\epsilon}}\right)$$
The reduced solution is 
$$u(x,y) = x$$
The solution has a single boundary layer at the outflow boundary $y=1$, as seen in Figure \ref{fig:2dsol}.

We fix the network structure to one hidden layer of width 5. We select residual samples as an equispaced mesh of size 129 by 129 over $[-1,1]\times [-1,1]$. Then we train the networks with different input transformation settings. It is evident that with the input transformations, the networks were able to capture the steep boundary layer at $y=1$ more effectively. 

\begin{figure}[t]
\centering
\begin{subfigure}{0.48\textwidth}
\centering
\includegraphics[width = \linewidth]{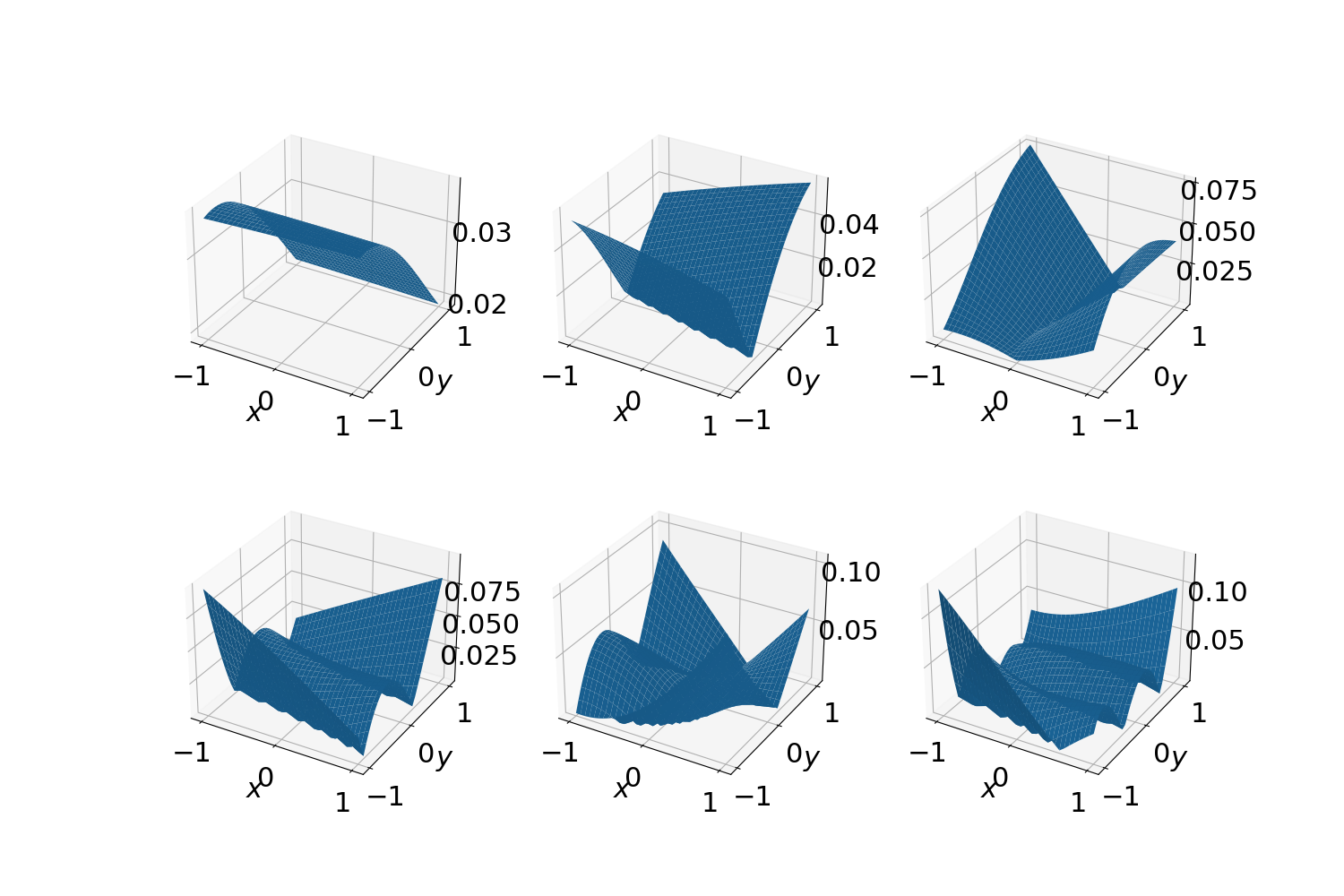}
\caption{}
\label{fig:2deigvec_none}
\end{subfigure}
\begin{subfigure}{0.48\textwidth}
\centering
\includegraphics[width = \linewidth]{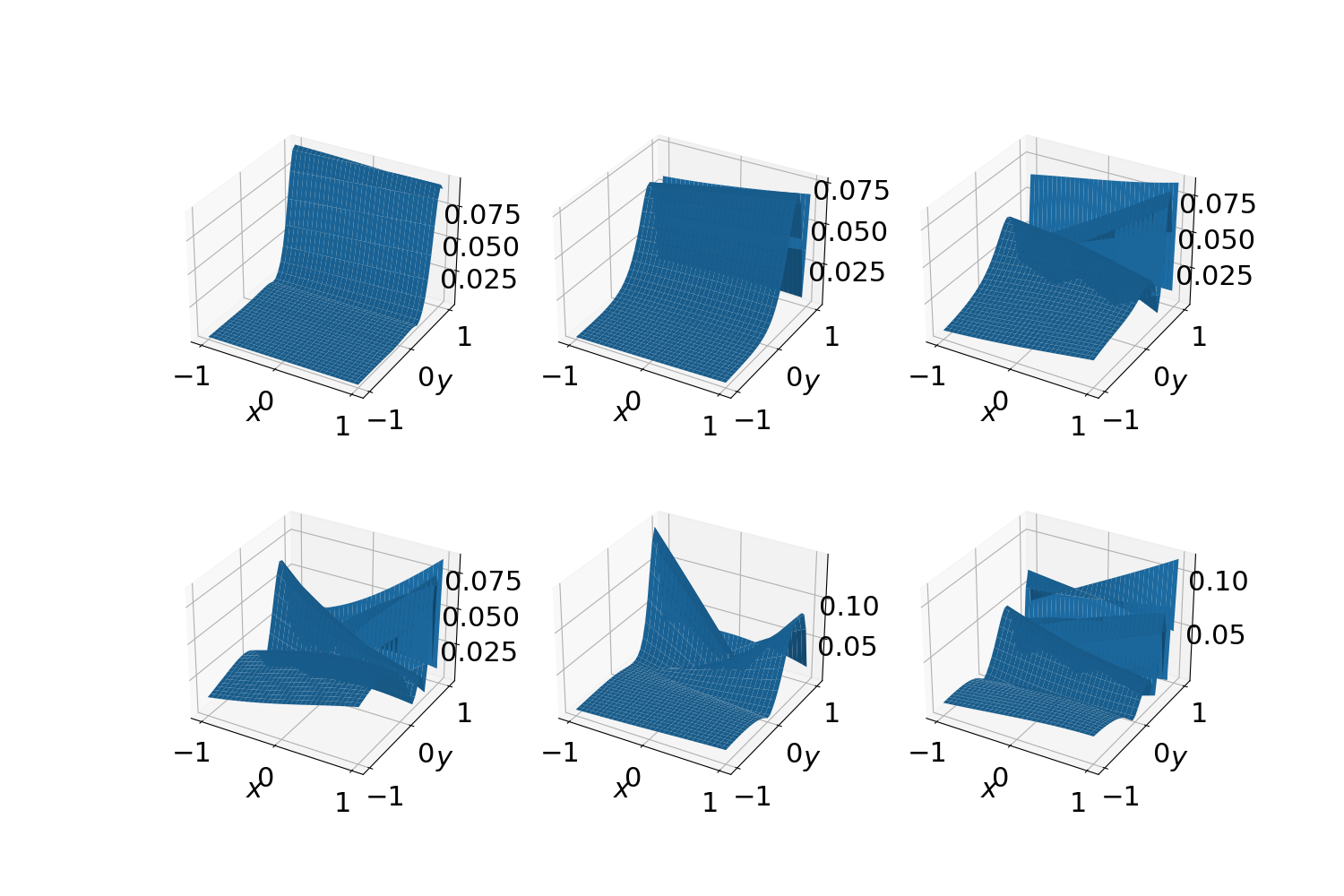}
\caption{}
\label{fig:2deigvec_inputm1}
\end{subfigure}
\caption{Eigenvectors corresponding to the 6 largest eigenvalues of $\bfK_{rr}$. Eigenvalues decrease from left to right, top to bottom. Fig \ref{fig:2deigvec_none}: the eigenvectors for no transformation. Fig \ref{fig:2deigvec_inputm1}: the eigenvectors for transformation $(x,y)\rightarrow (x,10(y-1))$. }
\label{fig:2deig_vec}
\end{figure}

We examine the behavior of eigenvectors of NTK of the two-dimensional problem to confirm the effect of input transformation. In terms of one-dimensional cases, when we perform input transformation $(x,y)\rightarrow (x,10(y-1))$, the first eigenvector forms a steep boundary layer near $y=1$. As a result, approximations by PINNs are more accurate with such input transformation. 

\section{Discussion}

This work examines the performance of standard physics-informed neural networks when applied to convection-diffusion equations, a representative problem within the class of singularly perturbed problems. Specifically, we observe how naturally PINNs tend to produce erroneous approximations to uncomplicated one-dimensional convection-diffusion equations. We offered a straightforward linear transformation technique that tunes the networks to accurately capture the behavior of solutions of convection-diffusion equations. The linear input transformation led the networks to arrive at the correct approximations, and it also reduced the number of training epochs needed. Trained PINNs employing linear input transformations compute solutions using a relatively small number of data points, for which traditional methods produce oscillatory non-physical solutions. We then examined the influence of this linear transformation technique and offered explanations through the lens of neural tangent kernels. This linear transformation technique does not need to be restricted to convection-diffusion equations. Other problems under the singularly perturbed problem class, such as reaction-diffusion equations, can also benefit from this technique. The ideas presented so far are focused on a one-dimensional convection-diffusion equation, but these concepts can be extended to two dimensions, and we have also demonstrated the application of linear input transformations on a two-dimensional example. 

% There are still many problems left unanswered. How does this linear transformation trick perform on two-dimensional convection-diffusion equations? Also, will the dynamics of the neural networks behave differently when the dimension of the problem increases? We still need to make strides in these directions to understand the behavior of PINNs for convection-diffusion equation. Comprehending the behavior of PINNs will lead us to designing better network structures for singularly perturbed equations. 

% \pagebreak
\printbibliography

% \bibliographystyle{siamplain}
% \bibliographystyle{unsrtnat}
% \bibliography{doc_ref} 

% \begin{thebibliography}{99}

% \bibitem{1} Spiegel, M. R. (1981). Theory and problems of Advanced Calculus: Si (metric) edition. McGraw-Hill. 

% \end{thebibliography}
\end{document}